\documentclass{amsart}
\usepackage{bbold}
\usepackage[arrow,matrix]{xy}
\usepackage[amsrefs]{math}

\newcommand\T{{\mathcal T}}
\newcommand\charin{{\lessdot}}
\newcommand\out{{\operatorname{\mathsf{Out}}}}
\newcommand\cent{{\operatorname{\mathsf{Cent}}}}
\newcommand\6{{\mathbb1}}
\newcommand\7{{\mathbb2}}
\newcommand\8{{\mathbb3}}
\newcommand\9{{\mathbb d}}
\renewcommand\sym[1]{{\mathsf{Sym}(#1)}}
\newcommand\DIAG{{\boldsymbol\bigtriangleup}}
\newcommand\PROD{{\boldsymbol\times}}
\newcommand\DERI[1][]{{\raise0.4ex\hbox{$#1\boldsymbol\bigtriangledown$}}}
\newcommand\PROJ{{\boldsymbol\dashv}}

\begin{document}
\author{Laurent Bartholdi}
\email{laurent.bartholdi@epfl.ch}
\address{\'Ecole Polytechnique F\'ed\'erale, SB/IGAT/MAD, B\^atiment
  BCH, 1015 Lausanne, Switzerland}
\author{Said N. Sidki}
\email{sidki@mat.unb.br}
\address{Universidade de Bras\'\i lia\\
  Departamento de Matem\'atica\\ 70.910-900 Brasilia-DF\\ Brasil} 
\title{The Automorphism Tower of Groups acting on Rooted Trees}
\thanks{The authors gratefully acknowledge support from the ``Fonds
  National Suisse de la Recherche Scientifique''.}
\date{August 15, 2003}
\begin{abstract}
  The group of isometries $\aut(\T_n)$ of a rooted $n$-ary tree, and
  many of its subgroups with branching structure, have groups of
  automorphisms induced by conjugation in $\aut(\T_n)$. This fact has
  stimulated the computation of the group of automorphisms of such
  well-known examples as the group $\Gg$ studied by R.~Grigorchuk, and
  the group $\GS$ studied by N.~Gupta and the second author.
  
  In this paper, we pursue the larger theme of towers of automorphisms
  of groups of tree isometries such as $\Gg$ and $\GS$. We describe
  this tower for all subgroups of $\aut(\T_2)$ which decompose as
  infinitely iterated wreath products. Furthermore, we describe fully
  the towers of $\Gg$ and $\GS$.
  
  More precisely, the tower of $\Gg$ is infinite countable, and the
  terms of the tower are $2$-groups. Quotients of successive terms are
  infinite elementary abelian $2$-groups.
  
  In contrast, the tower of $\GS$ has length $2$, and its terms are
  $\{2,3\}$-groups. We show that $\aut^2(\GS)/\aut(\GS)$ is an
  elementary abelian $3$-group of countably infinite rank, while
  $\aut^3(\GS)=\aut^2(\GS)$.
\end{abstract}
\maketitle

\section{Introduction}
The completeness of a centerless\footnote{in the sense that $G$ has a
  trivial center} group is measured by the extent to which all its
automorphisms are all inner, i.e.\ to which $G=\aut G$.  If $G$ has
outer automorphisms, then adjoining them to $G$ is part of the process
of completing $G$; indeed $G$ embeds in $\aut G$ as the subgroup of
inner automorphisms, and $\aut G$ is the one-step completion of $G$.
This process can be repeated on $\aut G$ to produce a subnormal series
of automorphism groups $\aut^iG$ for all $i\in\N$.  Formally, the
\emph{automorphism tower} of $G$ it is defined as follows:
\begin{defn}
  Given a centerless group $G$, we define $\aut^\alpha(G)$ for any
  ordinal $\alpha$ as follows:
  \begin{gather*}
    \aut^0(G)=G;\\
    \aut^{\alpha+1}(G)=\aut(\aut^\alpha(G));\\
    \intertext{and if $\alpha$ is a limit ordinal, then}
    \aut^\alpha(G)=\bigcup_{\beta<\alpha}\aut^\beta(G).
  \end{gather*}
  All these groups are also centerless, and hence form an ascending
  tower of groups.
  
  We define the \emph{height} $\tau(G)$ as the least ordinal $\alpha$
  such that $\aut^\alpha(G)=\aut^{\alpha+1}(G)$; this minimal $\alpha$
  exists by Theorem~\ref{thm:thomas}.

  The group $G$ is \emph{complete} if $\tau(G)=0$, i.e., if $\out(G)=1$.
\end{defn}

Many important groups have height one: for example, non-abelian free
groups $F$ satisfy $\aut^2 F=\aut F$; this was shown by J.~Dyer and
E.~Formanek~\cite{dyer-f:autf} if $F$ has finite rank, and by
V.~Tolstykh~\cite{tolstykh:at} for general $F$. Mapping class groups,
simple groups, and arithmetic subgroups of simple Lie groups, are also
of height one~\cite{bridson-v:outer}. The isometry group of a regular,
\emph{non-rooted} tree is complete~\cite{znoiko:automorphism}.

In this paper, we are more specifically interested in the automorphism
tower of groups acting on rooted, regular trees. These groups are
quite interesting in that their automorphism tower turns out often to
be quite short ($1$, $2$ or $\omega$, much less than the upper bound
the successor of $2^\omega$ given by Theorem~\ref{thm:thomas}), and
can in concrete cases be described explicitly as groups acting on the
original tree; that is, automorphisms of these groups are induced by
conjugation by tree isometries. Therefore, we deal with the normalizer
tower of a group defined as follows.
\begin{defn}  
  Given a subgroup inclusion $G\le N$, we define the \emph{normalizer
    tower} of $G$ in $N$ as
  \begin{gather*}
    \norm^0_N(G)=G;\\
    \norm^{\alpha+1}_N(G)=\norm_N(\norm_N^\alpha(G));\\
    \intertext{and if $\alpha$ is a limit ordinal, then}
    \norm^\alpha_N(G)=\bigcup_{\beta<\alpha}\norm^\beta_N(G).
  \end{gather*}
\end{defn}

To fix notation, let $\T$ be a $d$-regular rooted tree, for some
integer $d\ge2$. Let also $G$ be a group acting on $\T$. Consider the
$d$ subtrees $\T_1,\dots,\T_d$ rooted at all vertices neighbouring the
root of $\T$. The subtrees $\T_1,\dots,\T_d$ are permuted by $G$, and
the (set-wise) stabilizer of $\T_i$ acts on $\T_i$ by restriction.
Each of the subtrees $\T_i$ is isomorphic to $\T$, and therefore
carries an action of $G$ as well as a restricted action of
$\stab_{\T_i}(G)$ to the subtree $\T_i$. We call $G$ \emph{layered} if
the direct product of $d$ copies of $G$, each acting individually on
$\T_1,\dots,\T_d$, embeds in $G$, and $G$ permutes transitively the
subtrees $\T_1,\dots,\T_d$. A first goal of this paper is to give a
quite explicit way of calculating the automorphism tower of a layered
group (see Section~\ref{sec:layer}):
\begin{thm}\label{thm:mainL}
  Let $L$ be a layered group of isometries of the binary tree, set
  $N=\norm_{\aut\T}(L)$, and let $N_*$ be the \emph{layered closure}
  of $N$, i.e.\ the smallest layered group containing $N$. Then
  $\aut^i(L)\le N_*$ for all $i\le\omega$.
  
  More precisely, let $\mathcal C$ be the smallest
  lattice\footnote{i.e., family of groups closed under the operations
    $(H,K)\mapsto H\cap K$ and $(H,K)\mapsto\langle H,K\rangle$} of
  subgroups of $\aut\T$ containing $L$, $N$ and closed under the
  operations $G\mapsto G\times G$, $G\mapsto\setsuch{(g,g)}{g\in G}$,
  $G\mapsto\norm_N(G)$, and $(G,H)\mapsto\setsuch{g\in H}{g^2\in
    G\text{ and }[H,g]\le G}$ for every $G,H\in\mathcal C$ with
  $G\triangleleft H$. Then $\aut^i(L)\in\mathcal C$ for all $i\in\N$.
\end{thm}

Consider now a slightly weaker condition than being layered: that the
intersection of the groups $G$ and $G\times\dots\times G$ have finite
index in both groups. There are some finitely-generated examples of
such groups; the best-studied are the ``Grigorchuk group'' $\Gg$
acting on the binary tree~\cite{grigorchuk:burnside} and the group
$\GS$ acting on the ternary tree~\cite{gupta-s:burnside}.  Both of
these groups enjoy many extra properties, such as being torsion,
having intermediate word-growth, and being just-infinite. The
automorphism group of $\GS$ was described in~\cite{sidki:subgroups}.

A second goal of this paper is to describe quite explicitly all terms
in the automorphism tower of these two groups and groups similar to
them. The automorphism tower lies inside the full isometry group
of the original tree, and the last term of the tower has a natural
description:
\begin{thm}
  The group $\GS$ defined in~\cite{gupta-s:burnside} has height
  $2$.
\end{thm}

\begin{thm}
  The Grigorchuk group $\Gg$ has height $\omega$, and
  $\aut^\omega(\Gg)$ is the smallest layered subgroup of $\aut\T$
  containing $G$; equivalently, it is the subgroup of $\aut\T$
  generated by a copy of $G$ acting below each vertex.
\end{thm}
%

Finally, in a more general direction, we are interested in knowing
when the automorphism group of a countable group is again countable;
and how far the automorphism tower of such a group grows. An example
is the following result, whose proof closely
follows~\cite{thomas:atbook}*{Theorem~3.3.1}. We make freely use of
notions defined in that reference, in particular stationary sets
(Definition~3.3.4) and Fodor's lemma (Theorem~3.3.11). Roughly
speaking, a stationary set should be thought of as having positive
measure, and Fodor's lemma states that a strictly decreasing function
on a stationary set has a stationary fiber.
\begin{prop}\label{prop:towercountable}
  Let $H$ be a countable group. Suppose $G\le H\le K$ with $G$ finitely
  generated and $\cent_K(G)=1$. Then the terms of the normalizer tower
  of $H$ in $K$ are countable, and therefore the tower has countable
  height.
\end{prop}
\begin{proof}
  Assume $G=\langle s_1,\dots,s_k\rangle$, and denote the normalizer
  tower of $H$ by $\{N_\alpha\}$; i.e.\ $N_0=H$ and
  $N_{\alpha+1}=\norm_K(N_\alpha)$ and
  $N_\alpha=\bigcup_{\gamma<\alpha}N_\gamma$ for limit ordinals
  $\alpha$.
  
  First, we show that $N_\alpha$ is countable for all
  $\alpha<\omega_1$, the first uncountable ordinal. This is clearly
  true for $\alpha=0$. Then, for any $j\in\{1,\dots,k\}$, we have
  $s_j^{N_{\alpha+1}}\subseteq N_\alpha$ and therefore
  $|s_j^{N_{\alpha+1}}|\le\aleph_0$. This implies
  $[N_{\alpha+1}:\cent_{N_{\alpha+1}}(s_j)]\le\aleph_0$, and hence
  \[[N_{\alpha+1}:\bigcap_{i=1}^k\cent_{N_{\alpha+1}}(s_j)]=[N_{\alpha+1}:\cent_{N_{\alpha+1}}(G)]=|N_{\alpha+1}|\le\aleph_0.\]
  Finally, if $\alpha$ is a countable limit ordinal, then $N_\alpha$
  is a countable union of countable groups and is therefore countable.
  
  Assume now by contradiction that $N_\alpha\neq N_\beta$ for all
  distinct $\alpha,\beta<\omega_1$. Let $T$ be the set of countable
  limit ordinals, which is a stationary subset of $\omega_1$. For each
  $\alpha\in T$ choose an $h_\alpha\in N_{\alpha+1}\setminus
  N_\alpha$. Since $\alpha$ is a limit ordinal and
  $N_\alpha\triangleleft N_{\alpha+1}$, we have
  \[G^{h_\alpha}\le N_\alpha=\bigcup_{\gamma<\alpha}N_\gamma;\]
  now for each $j\in\{1,\dots,k\}$ choose $\gamma_j<\alpha$ such that
  $s_j^{h_\alpha}\in N_{\gamma_j}$, and set
  $f(\alpha)=\max\{\gamma_j\}$. We have $f(\alpha)<\alpha$ for all
  $\alpha\in T$, i.e.\ $f$ is a regressive function $T\to\omega_1$.
  By Fodor's lemma there exists $\gamma<\omega_1$ with
  \[S=\setsuch{\alpha\in T}{f(\alpha)=\gamma}\text{ stationary.}\]
  For each $\alpha\in S$ define $\phi_\alpha:G\to H_\gamma$ by
  $x\mapsto x^{h_\alpha}$. Any $\phi_\alpha$ is determined by the
  images of $G$'s generators, so
  \[|\setsuch{\phi_\alpha}{\alpha\in S}|\le|H_\gamma|^k<\omega_1;\]
  however $|S|=\omega_1$ because $S$ is stationary, so there are
  $\alpha_1\neq\alpha_2\in S$ with $\phi_{\alpha_1}=\phi_{\alpha_2}$;
  in other words, $1\neq h_{\alpha_1}h_{\alpha_2}^{-1}\in
  \cent_K(W)=1$ and we have reached a contradiction.
\end{proof}

Recall that a tree isometry $\alpha\in\aut\T$ is \emph{finite state}
if its action is induced by a finite transducer, and that the set of
finite-state isometries forms a group $R$ --- see
page~\pageref{def:finite-state} for the definition.
\begin{cor}
  Let $R$ be the group of finite state isometries. Then $\aut^i(R)$ is
  countable for all $i\in\N$, so $\tau(R)$ is countable.
\end{cor}
\begin{proof}
  Apply Proposition~\ref{prop:towercountable} with $K=\aut\T$ and $G$
  either the Grigorchuk group (if $d=2$) or the Gupta-Sidki group (if
  $d\ge3$). By Theorem~\ref{thm:ln}, we have
  $\aut^\alpha(R)=\norm_K^\alpha(R)$.
\end{proof}

As another illustration, consider $G\le\aut\T$ a finitely generated
group acting on the binary tree $\T=\{\6,\7\}^*$; and assume that $G$
contains for all $v\in\T$ an element $g_v$ fixing $v\6$ and $v\7$, and
having non-conjugate actions on the subtrees rooted at these vertices.
Then $\cent_{\aut\T}(G)=1$, and Proposition~\ref{prop:towercountable}
applies.

\section{The automorphism tower}
For a good survey of known results on automorphism groups,
see~\cite{romankov:automorphisms}. Even though the rest of the paper
does not rely on the considerations in this section, we include here
some relevant facts.

First, we note that the automorphism tower has been computed for some
classes of groups, and indeed that in many cases the tower has small
height.  W.~Burnside showed in~\cite{burnside:theory}, page~95
that $\aut(G)$ is complete if and only if $G$ is characteristic in
$\aut(G)$, i.e.\ if and only if $G\triangleleft\aut^2(G)$.
He then showed (op.\ cit., next page) that if $G$ is a non-cyclic
simple group, then $\tau(G)\le1$.  Furthermore, the symmetric groups
$\sym n$ have height $0$, when $n\ne6$.

H.~Wielandt showed in~\cite{wielandt:tower} that the tower of a
centerless, finite group is finite. Little was then known on the
height of infinite groups.

A.~Rae and J.~Roseblade showed in~\cite{rae-r:extremal} that if $G$ is
a \v Cernikov group, i.e.\ $G$ is a finite extension of an abelian
group and satisfies the minimal condition on subgroups, then $\tau(G)$
is finite. Then J.~Dyer and E.~Formanek showed in~\cite{dyer-f:autf}
that non-abelian free groups have height $1$; i.e.\ 
$\out(\aut(F_n))=1$ for $n>1$.

However, groups with infinite height abound among infinite groups; for
instance, the infinite dihedral group $D_\infty=\langle
a,b|a^2,b^2\rangle$ has an outer automorphism exchanging $a$ and $b$;
and $\aut D_\infty\cong D_\infty$, so $\tau(D_\infty)\ge\omega$.
J.~Hulse showed in~\cite{hulse:poly} that centerless polycyclic groups
have countable height (though not necessarily $\omega$). He showed that
$D_\infty$ has height $\omega+1$, by computing its automorphism tower.
It can be described using subgroups of $GL(2,\Z[\frac12])$ as follows:
\[\aut^\alpha(D_\infty)=\begin{cases}
  \begin{pmatrix}1&\frac1{2^\alpha}\Z\\0&\pm1\end{pmatrix} & \text{
    for }\alpha<\omega,\\
  \begin{pmatrix}1&\Z[\frac12]\\0&\pm1\end{pmatrix} & \text{
    for }\alpha=\omega,\\
  \begin{pmatrix}1&\Z[\frac12]\\0&\pm2^\Z\end{pmatrix} & \text{
    for }\alpha>\omega.\end{cases}
\]

The well-definedness of the height of a centerless group is a result
of S.~Thomas:
\begin{thm}[\cite{thomas:autt}]\label{thm:thomas}
  Let $G$ be a centerless group. Then there exists an ordinal
  $\alpha$, less than the successor of $2^{|G|}$, such that
  $\aut^\alpha(G)=\aut^{\alpha+1}(G)$.
  
  Furthermore, for any ordinal $\alpha$ there exists a centerless
  group $G$ with $\tau(G)=\alpha$.
\end{thm}
(This last statement was later improved by Just, Shelah and
Thomas~\cite{just-s-t:autt}.)

If $G$ has a non-trivial center, more care is needed, since the
$\aut^\alpha(G)$ no longer form a nested tower. In the definition of
$\aut^\alpha(G)$, one replaces the union by a directed limit.
J.~Hamkins proved in~\cite{hamkins:tower} that all groups $G$ have a
well-defined height, by actually showing that $\aut^\alpha(G)$ is
centerless for some ordinal $\alpha$.  Note that for groups possibly
with non-trivial center there is no explicit bound on the height of
the tower as a function of $|G|$, as in Theorem~\ref{thm:thomas}.

The considerations we make for the tower of automorphisms apply also
to Lie algebras, if ``automorphism group'' is replaced by ``algebra of
derivations''.

We cannot resist the temptation of repeating the main steps of the
proof of Theorem~\ref{thm:thomas}, since they show a strong relation
to groups acting on rooted trees:
\begin{proof}[Proof of Theorem~\ref{thm:thomas} (Sketch)]
  Define inductively the groups $W_\alpha$, for ordinals $\alpha$, by
  $W_0=\Z/2$; $W_{\alpha+1}=W_\alpha\wr(\Z/2)$, in which $W_\alpha$
  embeds as $W_\alpha\times1$; and for limit ordinals $\alpha$,
  $W_\alpha=\bigcup_{\beta<\alpha}W_\beta$.
  
  If $\alpha$ is of the form $1+\beta$, then the normalizer tower of
  $W_0$ in $W_\alpha$ stabilizes after exactly $\beta$ steps.
  
  By a result of Fried and Koll\'ar~\cite{fried-k:autfield}, there
  exists a (usually non-Galois) field $\mathbb K$ whose automorphism
  group is $W_\beta$.  Consider then the group $G=PGL_2(\mathbb
  K)\rtimes W_0$. The automorphism tower of $G$ parallels the
  normalizer tower of $W_0$ in $W_\beta$, in that
  $\aut^\gamma(G)=PGL_2(\mathbb K)\rtimes\norm^\gamma_{\aut\mathbb
    K}(W_0)$ for all $\gamma$.  It follows that the automorphism tower
  of $G$ stabilizes after $\alpha$ steps.
\end{proof}

There are still many open questions in the topic of group
automorphisms, and in particular, as to how large the automorphism
group can be.

Let us mention that there are finitely presented groups with
infinitely generated automorphism group, see for
instance~\cite{lewin:autig}, where the group
$(\Z[\frac12]\times\Z[\frac12])\rtimes\Z^2$ is shown to have that
property for an appropriate action of $\Z^2$; see
also~\cite{mccullough:autig}.  The groups $\Gg,\GS$ that we consider
in this paper also have the property of being finitely generated
torsion groups but having an infinitely generated outer automorphism
group; however, they are not finitely
presented~\cites{lysionok:pres,sidki:pres}.

\subsection{Structure of the paper}
In Section~\ref{sec:trees} we recall standard notation for groups
acting on rooted trees. We also recall results by Y.~Lavreniuk and
V.~Nekrashevych showing that, under certain conditions often satisfied
in practice, the terms of the automorphism tower act on the same tree
as the original group.


In Section~\ref{sec:layer}, we study in some detail the automorphism
tower of ``layered'' groups of isometries of the binary tree, and
compute these towers for certain concrete examples.

In Section~\ref{sec:gs}, we study the automorphism tower of the group
$\GS$ determined in~\cite{gupta-s:burnside}, and of its kin the groups
$\BG,\FG$ studied in~\cite{bartholdi-g:parabolic}. The last group
$\FG$ was first studied in~\cite{fabrykowski-g:growth1}.  The main
result is that, for the groups $\FG$ and $\GS$, the tower terminates
after two steps, i.e.\ they have height $2$.

Finally, in Section~\ref{sec:grig}, we study the automorphism tower of
the Grigorchuk group $\Gg$ acting on $\T=\{\6,\7\}^*$, and describe
the terms of the tower by their action on same tree $\T$. The last
term of the tower is shown to be the smallest layered group containing
$\Gg$.


\subsection{Notation} We say that $H$ is \emph{characteristic} in $G$,
written $H\charin G$, if $H^\phi=H$ for all $\phi\in\aut(G)$. If $G$
is centerless, we identify $G$ with $\mathsf{Inn}(G)$, and then
$H\charin G$ if and only if $H\triangleleft\aut(G)$. We then say that
$H$ is \emph{absolutely characteristic} in $G$ if
$H\triangleleft\aut^\alpha(G)$ for all ordinals $\alpha$.

We write $N:Q$ for a split extension of $N$ by $Q$, and $N\cdot Q$ for
an nonsplit extension.

As usual, we write $[a,b]$ for the commutator $a^{-1}b^{-1}ab$, and
$a^b=b^{-1}ab$ for the conjugation action.

\section{Groups acting on trees}\label{sec:trees}
Let $X=\{\6,\dots,\9\}$ be an alphabet with $d\ge 2$ elements. The
\emph{rooted tree} $\T$ on $X$ has as set of vertices the free monoid $X^*$,
namely the set of finite words $x_1\dots x_n$, with $x_i\in X$ for all
$i$. The tree structure is obtained by connecting $x_1\dots x_n$ to
$x_1\dots x_nx_{n+1}$ for all choices of $x_i\in X$; there is a
distinguished \emph{root vertex} $\emptyset$, the empty word. The
\emph{level} of the vertex $v=x_1\dots x_n$ is $|v|=n$. All vertices
of level $n$ form the $n$th \emph{layer} $X^n$ of $\T$.

The boundary $\partial\T$ of $\T$ is naturally the set $X^\omega$ of
infinite sequences over $X$. Given a \emph{ray}
$v=x_1x_2\dots\in\partial\T$, we denote by $v_n$ the truncation
$x_1\dots x_n$.


For $v\in\T$ of length $n$, we denote by $v\T$ the subtree of $\T$
spanned by all vertices $vw$ with $w\in\T$. Abstractly, it is a tree
rooted at $v$, isomorphic to $\T$.

Let $\sym X$ be the symmetric group on $X$. It acts naturally on $\T$
by
\begin{equation}\label{eq:rigid}
  (x_1x_2\dots x_n)^\sigma = x_1^\sigma x_2\dots x_n,
\end{equation}
and we will always identify $\sym X$ with its image in $\aut\T$.

Let $W$ denote the isomorphism group of $\T$. To avoid confusion with
the automorphism group of $W$ itself or of its subgroups, elements of
$W$ will be called \emph{isometries}. Given $g\in G$, there is
$\sigma_g\in\sym X$ such that $g\sigma_g^{-1}$ fixes the vertices of
the first level of $\T$; then restriction to subtrees $x\T$ for all
$x\in X$ gives \emph{states} $g@x\in W$, for all $x\in X$. Conversely,
states $g@x$ and $\sigma_g\in\sym X$ can be assembled to give $g\in
W$. We therefore have a wreath product structure
\[W=W\wr_X\sym X,\quad g\mapsto(g@\6,\dots,g@\9)\sigma_g.\]
This process can be iterated; for $n\in\N$, $v\in X^n$ and $g\in W$ we
write $g@v$ the state of $g$ at vertex $v$, and ${}^n\sigma_g$ the
permutation action of $g$ on $X^n$. In other words, $g@v$ denotes the
action $g$ does on the subtree $v\T$ \emph{before} vertex $v$ is moved
by ${}^n\sigma_g$.

There is a dual construction: for $v\in\T$ and $g\in W$ we write $v*g$
the isometry of $\T$ that acts as $g$ on $v\T$ and fixes all other
vertices. We then have the simple
\begin{lem}\label{lem:@*}
  $@$ is a right action of $X^*$ on $W$, and $*$ is a left action: for
  all $g,h\in W$ and $v,w\in X^*$ we have
  \begin{xalignat*}{2}
    (g@v)@w &= g@(vw), &       v*(w*g) &= (vw)*g,\\
    (gh)@v  &= (g@v)(h@v^g), & v*(gh)  &= (v*g)(v*h),\\
    g       &= (v*g)@v, &      g       &= \left(\prod_{v\in X^n}v*(g@v)\right)\,^n\sigma_g,
  \end{xalignat*}
  where the $v*(g@v)$ mutually commute when $v$ ranges over the $n$th
  layer $X^n$.
\end{lem}
We also define a variant of the $*$-action, by defining isometries
$Zg$ for $Z\in\{\PROJ,\DERI,\DIAG\}^*$:
\begin{defn}\label{defn:xaction}
  For $g\in\aut\T$, we define
  \[\DIAG g=(g,\dots,g),\quad\DERI g=(1,\dots,1,g,g^{-1}),\quad\PROJ g=(1,\dots,1,g)=\9*g,\]
  where a tuple $(g_1,\dots,g_n)$ designates an element $h\in\aut\T$
  fixing the first level of $\T$ and satisfying $h@\mathbb i=g_i$.
  
  For $H\le\aut\T$ and $Z\in\{\DIAG,\DERI,\PROJ\}$, we define
  $ZH$ as $\setsuch{Zh}{h\in H}$. We also let $\PROD H$ denote the
  subgroup $H\times\dots\times H$ generated by $\setsuch{\mathbb
    i*h}{i\in X, h\in H}$.
\end{defn}

We will be concerned with various classes of groups which emerged in
recent years as important subgroups of $W$. The \emph{vertex
  stabilizer} $\stab_G(v)$ is the subgroup of $G$ fixing $v\in\T$. We
make the following definitions: a subgroup $G$ of $W$ is
\begin{description}
\def\makelabel#1{\emph{#1},}
\item[level-transitive] if $G$ acts transitively on $X^n$ for all
  $n\in\N$;
\item[weakly recurrent] if it is level-transitive, and
  $\stab_G(x)@x<G$ for all $x\in X$;
\item[recurrent] if it is level-transitive, and $\stab_G(x)@x=G$ for
  all $x\in X$;
\item[saturated] if $G$ contains for all $n\in\N$ a characteristic
  subgroup $H_n$ fixing $X^n$ such that $H_n@v$ is level-transitive
  for all $v\in X^n$;
\item[weakly branch] if $G$ is level-transitive, and $(v*G)\cap G$ is
  non-trivial for all $v\in\T$;
\item[weakly regular branch] if $G$ is level-transitive, and has a
  non-trivial normal subgroup $K$ with $x*K<K$ for all $x\in X$;
\item[branch] if $G$ is level-transitive, and $\langle (v*G)\cap G:v\in
  X^n\rangle$ has finite index in $G$ for all $n\in\N$;
\item[regular branch] if $G$ is level-transitive, and has a
  finite-index normal subgroup $K$ with $x*K<K$ for all $x\in X$;
\item[pre-layered] if $x*G<G$ for all $x\in X$;
\item[layered] if $G$ is level-transitive and pre-layered.
\end{description}
The subgroup $\rist_G(v)=v*G\cap G$ is called the \emph{rigid stabilizer} of
$v\in\T$.  The \emph{layered closure} of $G$ is the subgroup
$G_*=\langle v*G:v\in\T\rangle$ of $W$. The \emph{closure} of $G$ is
its topological closure $\widehat G$ as a subgroup of $W$; one has
\[
\widehat G=\setsuch{g\in W}{\text{ for all }n\in\N\text{ there are
}g^{(n)}\in G\text{ with }g^{(n)}_n=g_n}.
\]
The topological closure of a layered group is always $W$.

The pointwise stabilizer in $G$ of the $n$th layer of $\T$ is written
$\stab_G(n)$. Clearly $G/\stab_G(n)$ is the permutation group on $X^n$
generated by $\setsuch{g_n}{g\in G}$.

Note that in a weakly branch group all rigid stabilizers are actually
infinite.

Clearly ``layered'' implies ``regular branch'', which implies
``branch'' and ``weakly regular branch'', each of which imply ``weakly
branch''.

The \emph{finitary group} $F<W$ is defined as the layered closure of
$\sym X$. It is a locally finite, countable subgroup of $W$, and hence
is a minimal layered subgroup. Every layered group contains $F$.

Fix an isometry $g\in W$. The \emph{activity} of $g$ at $v\in\T$ is
$\sigma_{g@v}\in\sym X$. The \emph{portrait} of $g$ is the activity map
$\T\to\sym X$, $v\mapsto \sigma_{g@v}$; there is a bijection between $W$
and the set of portraits.

\label{def:finite-state}
Consider finite-state automata with input alphabet $X$ and states
labelled by $\sym X$. The set of portraits defined by such automata
defines the group $R$ of \emph{finite-state isometries} of $\T$. Our
main examples of groups, and their automorphism groups, are subgroups
of $R$.

For a subgroup $A$ of $\sym X$, consider the subgroup $\widehat{A_*}$
of $W$ consisting of those elements whose activity at any vertex is in
$A$. If $A$ is transitive on $X$, the resulting group $\widehat{A_*}$
is layered. In particular, if $A=\langle(\6,\7,\dots,\9)\rangle$ and
$d$ is prime, then $C=\widehat{A_*}$ is a pro-$d$-Sylow of $W$. It
consists of all tree isometries whose action below any vertex belongs
to $A$.

\begin{lem}\label{lem:pre}
  Let $G$ be any group acting on $\T=X^*$. Then the set of pre-layered
  subgroups of $G$ forms a lattice of groups, closed under taking
  verbal subgroups.
\end{lem}
\begin{proof}
  If $H,K$ are pre-layered subgroups of $G$, then for any $g\in H\cap
  K$ we have $x*g\in H\cap K$, and for any $g=h_1k_1\dots h_nk_n\in
  \langle H,K\rangle$ we have $x*g\in\langle H,K\rangle$ by
  Lemma~\ref{lem:@*}.
  
  Finally, for any $W\in F(X_1,\dots,X_n)$ we set
  $W(H)=\setsuch{W(h_1,\dots,h_n)}{h_i\in H}$. Given
  $g=W(h_1,\dots,h_n)\in W(H)$, we have $x*g=W(x*h_1,\dots,x*h_n)\in
  W(H)$ so $W(H)$ is pre-layered.
\end{proof}

\begin{lem}
  If $L$ is a layered group acting on $\T=X^*$, then there exists a
  transitive subgroup $A$ of $\sym X$ such that $L\cong L\wr A$.
\end{lem}
\begin{proof}
  Let $A$ be the quotient of $L$ obtained by restricting the action of
  $L$ to $X\subset X^*$. Since $x*L<L$ for all $x\in X$, and
  $\stab_L(1)\le L^X$, we have $1\to\stab_L(1)=L^X\to L\to A\to 1$.
  Furthermore, every $a\in A$ has a lift $(\ell_1,\dots,\ell_d)a\in
  L$, and since $(\ell_1,\dots,\ell_d)\in L$ we have $(1,\dots,1)a\in
  L$, defining a splitting of $A$ in $L$.
\end{proof}

\begin{prop}[R. M\"oller~\cite{moller:automorphism}, Proposition~4]
  If $L$ is layered and $K\triangleleft L$ is a non-trivial normal
  subgroup, then $K$ contains $\stab_L(n)'$ for some $n\in\N$.
  
  Therefore every quotient of $L$ is abelian by finite.
\end{prop}

Assume now that $G$ is generated by a symmetric set $S$. This induces
a \emph{metric} on $G$, defined by
\begin{equation}\label{eq:metric}
  \|g\|=\min\setsuch{n}{g=s_1\dots s_n\text{ with }s_i\in S}.
\end{equation}
The group $G$ is \emdef{contracting} if there exist constants $\eta<1$
and $K$ such that $\|g@x\|\le\eta\|g\|+K$ for all $g\in\stab_G(1)$ and
$x\in X$.\footnote{A more general notion is usually
  used~\cite{bartholdi-g-n:fractal}, namely that there exists
  $\eta<1$, $n\in\N$ and $K$ such that $\|g@u\|\le\eta^n\|g\|+K$ for
  all $g\in G$ and $u\in X^n$. This notion becomes then independent of
  the choice of the generating set $S$, and there is furthermore an
  optimal $\eta$, which does not depend on $S$ either.}

\subsection{Automorphisms and isometries}
The following lemma is folklore:
\begin{lem}
  Let $G$ be a weakly branch group. Then the centralizer $\cent_W(G)$
  is trivial, and all conjugacy classes of $G$, except $\{1\}$, are
  infinite.
\end{lem}
\begin{proof}
  Let $a\in W$ be non-trivial; then $a$ moves a point $v\in\T$. Choose
  now any non-trivial $g\in\rist_G(v)$; then $[a,g]\neq1$, since it
  acts similarly to $g$ on $v\T$.
  
  Consider now a non-trivial $g\in G$. It moves a point $v\in\T$,
  hence its centralizer $\cent_G(g)$ intersects trivially
  $\rist_G(v)$, which is infinite, so $\cent_G(g)$ has infinite index
  in $G$ and $g^G$ is infinite.
\end{proof}

In~\cite{lavreniuk-n:rigidity}, Y.~Lavreniuk and V.~Nekrashevych prove
\begin{thm}\label{thm:ln}
  Let $G$ act on $\T$, and suppose $G$ is weakly branch. Then
  \[\norm_{\mathsf{Homeo}(\partial\T)}(G)=\aut(G).\]
  Assume moreover that $G$ is saturated. Then
  \[\norm_W(G)=\aut(G).\]
\end{thm}
\begin{proof}[Sketch of proof]
  The first step is to show that for any $\phi\in\aut(G)$ and $v\in\T$
  there is $w\in\T$ with $\rist_G(v)^\phi\ge \rist_G(w)$.
  
  Taking the shortest such $w$ induces a continuous map on the
  boundary $\partial\T$, which sends $v$ to $w$.
  
  If $G$ is saturated, then $|w|=|v|$, so this map is at each level of
  the tree a permutation of its vertices. Since $W$ is compact, there
  is a convergent subsequence of these permutations, which converges
  to a tree isometry.
\end{proof}
The proof is in fact essentially model-theoretic;
see~\cite{rubin:reconstruction}. Indeed it amounts to showing that the
algebraic structure of a saturated weakly branch group is sufficient
to reconstruct its action on $\T$.


\begin{lem}\label{lem:laysat}
  Let $G$ be a layered group acting on the binary tree. Then $G$ is
  saturated.
\end{lem}
\begin{proof}
  Define inductively $H_0=G$ and $H_{n+1}=\mho_2(H_n)$, the subgroup
  of $H_n$ generated by the squares of its elements. Then clearly
  $H_n$ is characteristic in $G$ for all $n$, and since $G$ is
  layered, $H_n$ contains for all $g\in G$ the element
  $\big(((\dots(g,1)^\sigma,1)^\sigma,\dots,1)^\sigma\big)^{2^n}=(g,\dots,g)$
  with $2^n$ copies of $g$; therefore $H_n$ fixes $X^n$ and acts
  transitively on every subtree at level $n$.
\end{proof}

\begin{lem}\label{lem:wbsat}
  Let $G$ be a group acting on the binary tree. If $G$ contains a
  weakly branch group, then $G$ is weakly branch.
  
  If $G$ contains a subgroup $H$ such that $(\mho_2)^n(H)$ acts
  transitively on every subtree at level $n$, for all $n\in\N$, then
  $G$ is saturated.
\end{lem}
\begin{proof}
  The first statement follows directly from the definition: if $G$
  contains the weakly branch subgroup $H$, then $(v*H)\cap H\neq1$ for
  all $v\in\T$, and hence $(v*G)\cap G\neq1$ for all $v\in\T$.
  
  For the second statement, $(\mho_2)^n(G)$ fixes $X^n$, and also acts
  transitively on every subtree at level $n$.
\end{proof}

\begin{prop}\label{prop:layeredsat}
  The groups $F$ and $W$ are layered and saturated.
\end{prop}
\begin{proof}
  Branchness is clear, since we have $\rist_W(v)=v*W$ and
  $\rist_F(v)=v*F$.
  
  To prove saturatedness, we invoke Theorem~\ref{thm:swreath}; indeed
  choose $x\in X$ and write $W=\rist_W(x)\wr_X\sym X$; then
  $H_1\charin W$ may be chosen to be its base group $W^X$.
  Inductively, let $K$ be the copy of $H_{n-1}$ inside
  $\rist_W(x)=x*W$; then $H_n\charin W$ may be chosen as $K^X$,
  characteristic in $W^X$ which is itself characteristic in $W$.

  The same argument applies to $F$.
\end{proof}

Wreath products exhibit a form of rigidity which was already noted
independently by a number of authors:
\begin{thm}[Peter Neumann, \cite{neumann:swreath}; Yurii Bodnarchuk,
  \cite{bodnarchuk:nswreath}; Paul Lentoudis and Jacques Tits, \cite{lentoudis-t:autwreath}]
  \label{thm:swreath}
  If $A,B$ are non-trivial groups with $B$ acting on a set $X$, then
  the base group $A^X$ is characteristic in the wreath product $A\wr_X
  B$, unless $A=\Z/2$ and $B$ has an abelian subgroup $B_0$ of index
  $2$ containing unique square roots of its elements.
\end{thm}

Furthermore, in some specific examples, a strengthening of this
rigidity has been obtained: in essence, not only are all group
isomorphisms induced by a tree isomorphism, but moreover there is a
unique minimal tree carrying the group's action. R.~Grigorchuk and
J.~S.~Wilson showed that if $G$ is a branch group acting on a group
$\T$ satisfying two technical conditions, and if $G$ also acts on
another tree $\T'$ as a branch group, then the actions on $\T$ and
$\T'$ are intertwined by map $\T'\to\T$ defined by the erasing of some
levels in $\T'$.

In contrast, finite groups behave in a manner very different from the
groups we are interested in:
\begin{thm}
  If the finite group $G$ acts on $\T$, then $\cent_W(G)$ is uncountable.
\end{thm}
\begin{proof}
  For simplicity we prove the result only for $\T$ the binary tree
  $\{\6,\7\}^*$. We proceed by induction, first on $|G|$, and then on
  the smallest level at which a non-trivial activity occurs, i.e.\ the
  first $k$ such that $G$ is not a subset of $\stab_W(k)$.
  
  If $|G|=1$, then $\cent_W(G)=W$ and we are done. If
  $G\le\stab_W(1)$, then consider the projections $H_\6,H_\7$ of
  $G$ defined by restriction to the respective subtrees rooted at
  level $1$.  These are finite groups, of size at most $|G|$, and at
  least one of $H_\6,H_\7$, say $H_\6$, has already been covered by
  induction; therefore $\cent_W(G)$ contains
  $\cent_W(H_\6)\times\{1\}$ and is uncountable.
  
  We may therefore assume that $G$ contains an element
  $g=(g_\6,g_\7)\sigma$; let also $H$ be the projection of
  $\stab_G(1)$ on the subtree rooted at $\6$. Set $Z=\cent_W(H)$; by
  induction, $Z$ is uncountable. Furthermore, $\stab_G(1)\le H\times
  H^{g_\7^{-1}}$, so the centralizer of $\stab_G(1)$ contains $Z\times
  Z^{g_\7^{-1}}$.
  
  The centralizer of $g$ contains all elements $(z,z^{g_\7^{-1}})$ with
  $z\in\cent_W(g_\6g_\7)$; and since
  $g^2=(g_\6g_\7,g_\7g_\6)\in\stab_G(1)$, we see that $\cent_W(G)$
  contains all $(z,z^{g_\7^{-1}})$ for $z\in Z$; therefore
  $\cent_W(G)$ is uncountable.
\end{proof}

Consider also the following example of a group with uncountable
normalizer: let $U$ be an infinite minimal connecting set of vertices
of the tree $\T$, i.e.\ a subset of $\T$'s vertices that intersects
any infinite ray. For every $u\in U$, let $H(u)$ be a subgroup of
$\aut(\T)$. Let $G$ be the group generated by all $u*H(u)$, with $u$
ranging over $U$. Then $G$ is isomorphic to the direct sum of the
$H(u)$'s, and the normalizer of $G$ in $W$ contains the cartesian
product of the $H(u)$'s.

Let us note finally a result related to Theorem~\ref{thm:ln}, and
expressed in terms of commensurators: an \emph{almost automorphism} of
$G$ is an automorphism between two finite-index subgroups of $G$. Two
almost automorphisms are equivalent if they agree on on a finite-index
subgroup of $G$. The set of equivalence classes of almost
automorphisms carries a natural group structure, and is called the
\emph{abstract commensurator} of $G$.

On the other hand, if $G$ is a subgroup of $K$, then the commensurator
of $G$ in $K$ is $\setsuch{k\in K}{[G:G\cap G^k]<\infty,[G^k:G\cap
  G^k]<\infty}$.  
\begin{thm}[Claas R\"over, \cite{rover:commensurator}]
  The abstract commensurator of a weakly branch group $G$ acting on
  $T$ is isomorphic to the commensurator of $G$ in
  $\mathsf{Homeo}(\partial\T)$.
\end{thm}

\subsection{The automorphism tower}\label{subs:aut}
Our first purpose is, given a suitable $G$, to identify a group $G_*$
that is much smaller than $W$ (i.e., for instance, is countable); but
still is large enough to contain the whole automorphism tower of $G$.
This could be interpreted as a strengthening of Theorem~\ref{thm:ln}.

For any $n\in\N$ set $G_n=\langle v*G:\,v\in \T,|v|\le n\rangle$, and
$G_*=\bigcup_{n\ge0}G_n$. The following is immediate:
\begin{lem}
  Assume $G$ is weakly branch and recurrent. Then $G_n$ is weakly
  branch for all $n$, and the $G_n$ form an ascending tower. Its limit
  $G_*$ is layered, and is the smallest layered group containing $G$.
  If $G$ is countable, then so is $G_*$.
\end{lem}

\begin{thm}\label{thm:stops@omega}
  If $G$ is a finitely generated group with a recurrent, saturated,
  and weakly branch action on $\T$, and if $G$ satisfies
  $G_*=\bigcup_{n\in\N}\norm^n(G)$, then $G_*$ is complete; in other
  words, $G$ has height at most $\omega$.
\end{thm}
We remark that the condition $G_*=\bigcup_{n\in\N}\norm^n(G)$ implies that
for all $n\in\N$ there are $k,\ell\in\N$ such that
$\PROD^kG<\norm^n(G)<(\PROD^\ell G)(\wr^\ell C_2)$, since $G$ is finitely
generated.
\begin{proof}
  Take $\phi\in\aut(G_*)$. Then by Theorem~\ref{thm:ln} we have
  $\phi\in\aut\T$, and similarly $\psi=\phi^{-1}\in\aut\T$.  Since by
  definition $G$ is finitely generated, there exists $n\in\N$ such
  that $\phi(G)\le\norm^n(G)$ and $\psi(G)\le\norm^n(G)$.
  
  Since $G_*=\bigcup_{m\in\N}G_m$, there is also an $m\in\N$ such that
  $\phi(G)\le G_m$ and $\psi(G)\le G_m$.
  
  Let $u,\dots,v$ be all the vertices on level $m$ of $\T$, and let us
  decompose $\phi=(\phi_u,\dots,\phi_v)\sigma$ and
  $\psi=(\psi_u,\dots,\psi_v)\sigma^{-1}$. Since $G$ is recurrent, for
  any $g\in G$ there exists a $g'\in G$ with level-$m$ decomposition
  $(g,*,\dots,*)$. Then $(g')^\phi=(g^{\phi_u},*,\dots,*)^\sigma$ belongs to
  $G$, and therefore $\phi_u$ induces by conjugation an endomorphism
  of $G$.

  Then $(g,*,\dots,*)=g'=\phi\psi(g')=(\phi_u(*),*,\dots,*)$, and
  therefore $\phi_u(G)=G$, so $\phi_u$ induces an automorphism of
  $G$. We have shown $\aut(G_*)\le(\aut(G)F)_*=G_*$.
\end{proof}


%
%

\section{Layered groups}
\label{sec:layer}
We describe in this section the automorphism tower of a layered group
$L$ acting on the binary tree $\T=X^*$, with $X=\{\6,\7\}$; recall
that such group has a decomposition $L=(L\times L)\rtimes\langle
\sigma\rangle$, where $\sigma$ denotes the transposition of the top
two branches of the binary tree.

Denote $\aut(L)$ by $A$. By Lemma~\ref{lem:laysat} and
Theorem~\ref{thm:ln}, every automorphism $\alpha$ of $L$ acts on the
binary tree, and therefore we have $A\le(W\times
W)\rtimes\langle\sigma\rangle$. Then, given
$\alpha=(\alpha_\6,\alpha_\7)\sigma^i\in A$, we have
$(\ell,1)^\alpha=(\ell^{\alpha_\6},1)^{\sigma^i}\in L$, and similarly
$(1,\ell)^\alpha\in L$ so $\alpha_i\in A$ for all $i\in X$ and
$A\le(A\times A)\rtimes\langle\sigma\rangle$. Furthermore,
$\sigma^\alpha=\sigma^{1-i}(\alpha_\7^{-1}\alpha\6,\alpha_\6^{-1}\alpha_\7)\sigma^i\in
L$, so $\alpha_\7^{-1}\alpha_\6\in L$ and $A\le(A\times 1)L$.  In other
words, $\alpha$ can be written as a (possibly infinite) product
\[\alpha=\cdots(\6^n*\ell_n)\cdots(\6*\ell_1)\ell_0,\]
for some $\ell_i\in L$. On the other hand, such a product defines an
endomorphism (but not necessarily an automorphism) of $L$, by
conjugation.

Note also that $\alpha$ is finite state provided the sequence
$(\ell_n)$ is ultimately periodic (see~\cite{brunner-s:auto}). That
same paper shows that, for $F$ the finitary group, $\aut(F)$ contains
a copy of $W$ and is therefore uncountable.

We mention in passing another proof that $\aut(F)$ is uncountable:
recall that an automorphism $\phi$ of a group $G$ is \emph{locally
  inner} if $\phi_{|X}$ is inner for every finite $X\subset G$. Let
$\mathsf{Linn}(G)$ denote the group of locally inner automorphisms of
$G$.
\begin{prop}[\cite{hartley:slfg}, page 37]
  Let $G$ be a countably infinite, locally finite group, and suppose
  that $\cent_G(F)$ is not central in $G$, for all finite subgroups
  $F$ of $G$.  Then $|\mathsf{Linn}(G)|=2^{\aleph_0}$.
\end{prop}
It is, however, unknown whether $|\aut(G)|=2^{\aleph_0}$ for all
countably infinite, locally finite groups $G$.

\begin{defn}\label{defn:*}
  For a layered group $L$, we write $\aut^i(L)$ for the $i$-th term in
  the automorphism tower of $L$.  For groups $K\triangleleft G$, we
  write
  \[\Omega(G,K)=\{v\in G|\,v^2\in K\text{ and }[G,v]\leq K\}.\]
\end{defn}

We shall be concerned, in this section, with groups $G\le W$ that admit a
decomposition
\begin{equation}\label{eq:decomp}
  G=(\PROD U)\langle\sigma\rangle(\DIAG S)
\end{equation}
where $U\triangleleft S$ is maximal, and $G$ contains a fixed layered
group $L$. By Theorem~\ref{thm:ln} we have $\aut G=\norm_W(G)=:\norm
G$.  The form~\eqref{eq:decomp} is preserved by some important
operations, namely if $G=(\PROD U)\langle\sigma\rangle(\DIAG S)$ and
$H=(\PROD V)\langle\sigma\rangle(\DIAG T)$, then
\begin{align}
  G\cap H&=(\PROD(U\cap V))\langle\sigma\rangle(\DIAG(S\cap T)),\notag\\
  \Omega(G,H)&=(\PROD(U\cap\Omega(S,V)))\langle\sigma\rangle(\DIAG(\Omega(S,T))),\tag{Lemma~\ref{lem:WisY}}\\
  \norm(G)&=(\PROD\Omega(S,U))\langle\sigma\rangle(\DIAG(\norm
  U\cap\norm S)).\tag{Theorem~\ref{thm:normL}}
\end{align}

\begin{lem}\label{lem:norm}
  Consider a group $H=(\PROD K)\langle\sigma\rangle$. Then
  \[\norm(H)=(\PROD K)\langle\sigma\rangle(\DIAG\norm(K)),\]
  where $\PROD K$ is the largest geometrically decomposable subgroup
  of $\norm(H)$ and $\DIAG\norm(K)$ is the largest diagonal subgroup
  of $\norm(H)$.
\end{lem}
\begin{proof}
  Choose $\alpha=(\alpha_\6,\alpha_\7)\sigma^i\in \norm(H)$.  Then
  $(\alpha_\6,\alpha_\7)\in \norm(H)$, and $\alpha_\6,\alpha_\7\in \norm(K)$.
  Thus, $\alpha\in \PROD\norm(K)\langle\sigma\rangle$ and
  $[(\alpha_\6,\alpha_\7),\sigma]=(\alpha_\6^{-1}\alpha_\7,\alpha_\7^{-1}\alpha_\6)
  \in \PROD K$.  Therefore, $\alpha _2=\alpha_\6k$ for some $k\in K$; that
  is, $\alpha=(\alpha_\6,\alpha_\6k)\sigma^i$. Since
  $\DIAG\alpha_\6\in \norm(H)$ for all $\alpha_\6\in \norm(K)$, the
  decomposition $\norm(H)=(\PROD K)\langle\sigma\rangle(\DIAG\norm(K))$ follows.
  
  The rest of the assertions are easy.
\end{proof}

\begin{thm}\label{thm:normL}
  Consider $R=(\PROD U) \langle\sigma\rangle (\DIAG S)$ where $U$ is a
  normal subgroup of $S$. Then
  \[\norm(R)=(\PROD V) \langle\sigma\rangle (\DIAG T)\]
  where $V=\Omega(S,U)$ and $T=\norm(U)\cap \norm(S)$.
  
  Furthermore, $V$ is a normal subgroup of $T$, and $[\PROD V,R] \leq
  (\PROD U)(\DIAG V)$, and the quotient group $V/U$ is a $T/S$-module,
  and
  \begin{equation*}
    \norm(R)/R\cong (V/U) (T/S).
  \end{equation*}
\end{thm}
\begin{proof}
  First, it is clear that $V$ is a normal subgroup of $T$ and that
  $\PROD V$ normalizes both $\PROD U$ and $\DIAG S$.
  
  Choose $(v,1)\in \PROD V$; then
  $[(v,1),\sigma]=(v^{-1},v)=(v^{-2},1)( v,v)$ where by the definition
  of $V$ we have $v^{-2}\in U$ and $(v,v)\in\DIAG S$. Therefore $\PROD
  V\leq \norm(R)$. Easily, $\DIAG T\leq \norm(R)$.
  
  Choose next $\beta\in \norm(R)$. Then we may assume
  $\beta=(\beta_\6,\beta_\7)$. As $\PROD U$ is normal in $R$, we conclude
  that $\beta_\6,\beta_\7\in \norm(U)$ and commutation with $\DIAG S$
  shows that $\beta_\6,\beta_\7\in \norm(US)=\norm(S)$. Thus,
  $\beta_\6,\beta_\7\in T$ and
  \[\norm(R)\leq (\PROD T)\langle\sigma\rangle.\]
  
  Now, $\sigma^\beta=(\beta_\6^{-1}\beta_\7,\beta_\7^{-1}\beta_\6)
  \sigma=(u_\6s,u_\7s)\sigma$ for some $u_\6,u_\7\in U,s\in S$. Since
  $\sigma^2=1$, we have $u_\7s=(u_\6s)^{-1}=u_\6^{-s}s^{-1}$, and
  therefore
  \[s^2=u_\7^{-1}u_\6^{-s}\in U,\qquad \beta_\7=\beta_\6u_\6s.\]
  Thus, $\beta=(\DIAG\beta_\6)(1,u_\6s)$ with $\beta_\6\in T$.
  Since we have already shown $\DIAG T\leq \norm(R)$, and as $1\times
  U\leq R$, we may assume $\beta=(1,u_\6s)$, and as $1\times U\leq R$
  we may assume further that $\beta=(1,s)$. Now,
  $[\DIAG S,\beta]=1\times[S,s] \leq R$ and therefore $[S,s]\leq U$.
  Hence, $s\in\Omega(S,U)=V$ and
  $\norm(R)=(\PROD V)\langle\sigma\rangle(\DIAG T)$.

  We then have
  \begin{align*}
    [\PROD V,R] &= [\PROD V,\PROD U] [\PROD V,\DIAG S\langle\sigma\rangle]\\
    &=(\PROD [V,U]) (\PROD [V,S]) \setsuch{(v^{-1},v)}{v\in V}\\
    & \leq (\PROD U)(\DIAG V),
  \end{align*}
  since $[V,S],(V^2,1)\leq U$. We observe the following facts:
  \begin{gather*}
    (\DIAG T)\cap R=(\DIAG T)\cap (\DIAG(US))=\DIAG(T\cap S)=\DIAG S,\\
    ((\PROD V)\langle\sigma\rangle)\cap R=(\PROD U)(\DIAG V)\langle\sigma\rangle,\\
    V/U \text{ is a $T$-module and the kernel of the action contains $S$.}
  \end{gather*}
  
  Finally, let $W$ be a transversal of $U$ in $V$ and $X$ a
  transversal of $S$ in $T$. Then $(W\times 1) (\DIAG X)$ is a
  transversal of $R$ in $\norm(R)$, and hence $\norm(R)/R\cong (V/U)(T/S)$.
\end{proof}

\begin{rem}
  Suppose $\norm(S)=S$.  Then
  \[\norm^i(R)=(\PROD U_{i}) \langle\sigma\rangle(\DIAG S)\]
  where $U_0=U$ and $U_{i}=\Omega(S,U_{i-1})$.  The subgroups
  $U_{i}$ form a hypercentral series of $S$ with respect to $U$.
  
  In particular, if $S$ is finitely generated, then
  $\norm^\omega(R)=\norm^{\omega+1}(R)$, because letting
  $U_\omega=\bigcup_{n\in\N}U_n$ we have
  $\Omega(S,U_\omega)=U_\omega$.
\end{rem}

\begin{lem}\label{lem:WisY}
  Consider $G=(\PROD U)\langle\sigma\rangle(\DIAG S)$, and let $H=(\PROD V)
  \langle\sigma\rangle (\DIAG T)$ be a normal subgroup of $G$. Then
  \[\Omega(G,H)=(\PROD(U\cap\Omega(S,V))\langle\sigma\rangle (\DIAG\Omega(S,T)).\]
\end{lem}
\begin{proof}
  By assumption we have $V\triangleleft S$, $T\triangleleft S$,
  $U\triangleleft S$, $V\triangleleft T$, and $U\le\Omega(T,V)$.
  
  Set $Q=\Omega(G,H)$. First, consider $y=(x,1)$ with $x\in
  U\cap\Omega(S,V)$. Then $x^2\in V$ and $[x,S]\subseteq V$ so $y^2\in
  H$ and $[y,G]\in H$, so $y\in Q$. Conversely, if $y=(1,x)\in Q$,
  then $[y,G]\in H$, $[y,\sigma]\in H$ and $y^2\in H$ imply $x\in
  U\cap\Omega(S,V)$.
  
  Consider next $y=(x,x)$ with $x\in\Omega(S,T)$. Then $x^2\in T$ and
  $[x,S]\subseteq T$ so $y^2\in H$ and $[y,G]\in H$, so $y\in
  Q$. Conversely, if $y=(x,x)\in Q$, then $[y,G]\in H$ and $y^2\in H$
  imply $x\in\Omega(S,T)$.
\end{proof}

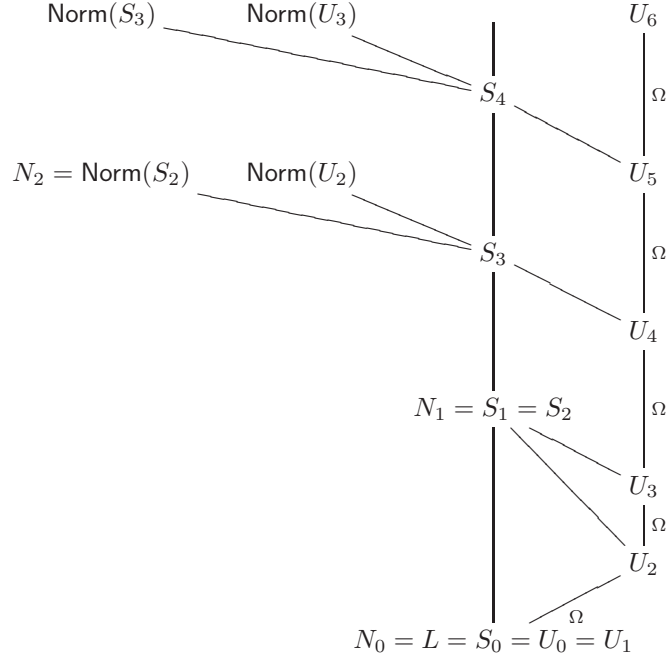
\begin{figure}
  \[\xymatrix@=15pt{
    {\norm(S_3)}\ar@{-}[drr] & {\norm(U_3)}\ar@{-}[dr] & \ar@{-}[d] & {U_6}\ar@{-}[dd]^{\Omega}\\
    & & {S_4}\ar@{-}[dd]\ar@{-}[dr]\\
    {N_2=\norm(S_2)}\ar@{-}[drr] & {\norm(U_2)}\ar@{-}[dr] & & {U_5}\ar@{-}[dd]^{\Omega}\\
    & & {S_3}\ar@{-}[dd]\ar@{-}[dr]\\
    & & & {U_4}\ar@{-}[dd]^{\Omega}\\
    & & {N_1=S_1=S_2}\ar@{-}[ddd]\ar@{-}[ddr]\ar@{-}[dr]\\
    & & & {U_3}\ar@{-}[d]^{\Omega}\\
    & & & {U_2}\ar@{-}[dl]^{\Omega}\\
    & & {\kern-3em N_0=L=S_0=U_0=U_1\kern-3em}
  }\]
  \caption{Some groups appearing in Theorem~\ref{thm:layered}}\label{fig:norms}
\end{figure}

\begin{thm}\label{thm:layered}
  Let $L$ be a layered group, $L=(\PROD L)\langle\sigma\rangle$ and let
  $(N_i=\norm^i(L))_{i\ge0}$ be the normalizer tower of $L$. Then we
  have the equations (see Figure~\ref{fig:norms})
  \begin{align*}
    N_0 &=U_0=S_0=L,\\
    N_i &=(\PROD U_i) \langle\sigma\rangle (\DIAG S_i),\\
    U_i &=\Omega(S_{i-1},U_{i-1}),\\
    S_i &=\norm(U_{i-1}) \cap \norm(S_{i-1}).
  \end{align*}
  
  For $1\leq i\leq 4$, the description of $N_i$ by the pair
  $(U_i,S_i) $ is as follows:
  \begin{xalignat*}{2}
    U_1 &=(\PROD L)\langle\sigma\rangle=L,          & S_1&=N_1;\\
    U_2 &=(\PROD L)\langle\sigma\rangle(\DIAG U_2), & S_2&=N_1;\\
    U_3 &=(\PROD L)\langle\sigma\rangle(\DIAG U_3), &
    S_3 &=(\PROD U_2)\langle\sigma\rangle(\DIAG N_1)=N_2;\\
    U_4 &=(\PROD U_2)\langle\sigma\rangle(\DIAG\Omega(N_1,U_3)), &
    S_4 &=(\PROD\Omega(U_3,U_1))\langle\sigma\rangle(\DIAG N_1).
  \end{xalignat*}
  Furthermore, the outer-normalizer groups are
  \begin{align*}
    N_2/N_1 &\cong U_2/U_1,\\
    N_3/N_2 &\cong (U_3/U_2) (N_2/N_1),\\
    N_4/N_3 &\cong (U_4/U_3) (\Omega(U_3,L)/U_2) (N_2/N_1),
  \end{align*}
  where $U_{i+1}/U_i$ and $\Omega(U_3,L)/U_2$ are elementary abelian
  $2$-groups.
\end{thm}
\begin{proof}
  Since $L=(\PROD L) \langle\sigma\rangle (\DIAG L)$, the general form
  of the normalizer subgroups $\norm^i(L)=(\PROD
  U_i)\langle\sigma\rangle(\DIAG S_i)$ is determined by
  Theorem~\ref{thm:normL}. We compute for $i\le4$:
  \subsection*{$\mathbf{i=1}$} We have $N_1=(\PROD U_1)\langle\sigma\rangle
  (\DIAG S_1)$, with $U_1=\Omega(L,L)=L=U_0$ and $S_1=\norm(L)\cap
  \norm(L)=\norm(L)=N_1$. We do not obtain any information about the
  outer-normalizer group $N_1/N_0$, since
  \[N_1/N_0\cong (U_1/U_0) (S_1/S_0) =S_1/S_0\cong N_1/N_0.\]
  \subsection*{$\mathbf{i=2}$} First, $N_2=(\PROD U_2) \langle\sigma\rangle
  (\DIAG S_2)$, with $U_2=\Omega (N_1,U_1)$ and $S_2=\norm(U_1)\cap
  \norm(S_1)$.

  Secondly, $S_2=\norm(L)\cap \norm(N_1)=N_1\cap N_2=N_1$.
    
  Thirdly, $N_0=L=L(\DIAG L)\leq U_2=\Omega(N_1,U_1)\leq N_1=L(\DIAG
  N_1)$. By Lemma~\ref{lem:WisY}, $U_2=L(\DIAG W)$, with $W=\Omega
  (N_1,N_0)=\Omega(N_1,U_1)=U_2$; therefore,
  \[U_2=L(\DIAG U_2)=U_1(\DIAG U_2).\]
    
  Fourthly, $\norm(U_2)=(\PROD V)\langle\sigma\rangle (\DIAG Y)$, with
  $V=\Omega(U_2,U_1)=U_2$ and $Y=\norm(U_1)\cap \norm(U_2)=N_1\cap
  \norm(U_2)=N_1$, and therefore $N_2=\norm(U_2)$.
    
  Finally, $N_2/N_1\cong (U_2/U_1) (S_2/S_1)=U_2/U_1$, an elementary
  abelian $2$-group.
  
  \subsection*{$\mathbf{i=3}$} First, $N_3=(\PROD U_3)
  \langle\sigma\rangle (\DIAG S_3)$, with $U_3=\Omega (S_2,U_2)$ and
  $S_3=\norm(U_2)\cap \norm(S_2)$.
    
  Secondly, $S_3=N_2\cap \norm(S_2)=N_2$.
    
  Thirdly, $U_3=\Omega(N_1,U_2)$ and $L(\DIAG U_2)\leq U_3\leq
  N_1=L(\DIAG N_1)$, so $U_3=L(\DIAG W)$ where
  $W=\Omega(N_1,U_2)=U_3$; and hence, $U_3=L(\DIAG U_3)$.
    
  Finally, $N_3/N_2\cong (U_3/U_2)(S_3/S_2)=(U_3/U_2)(N_2/N_1)$ where,
  as we had shown, $N_2/N_1\cong U_2/U_1$.

  \subsection*{$\mathbf{i=4}$} First, $N_4=(\PROD U_4) \langle\sigma\rangle
  (\DIAG S_4)$, with $U_4=\Omega (S_3,U_3)=\Omega (N_2,U_3)$ and
  $S_4=\norm(U_3)\cap \norm(S_3)$.

  Secondly, $S_4=\norm(U_3)\cap \norm(N_2)=\norm(U_3)\cap N_3$.
    
  Thirdly, from $U_3=(\PROD L) \langle\sigma\rangle (\DIAG U_3)$ we
  find
  \[\norm(U_3)=(\PROD W)\langle\sigma\rangle(\DIAG(\norm(L)\cap\norm(U_3))),\]
  with $W=\Omega(U_3,L)$.  Then $\norm(L)\cap \norm(U_3)=N_2\cap
  \norm(U_3) =S_3\cap \norm(U_3)=N_2$, and $\norm(U_3)=(\PROD W)
  \langle\sigma\rangle (\DIAG N_2)$. Since $W=\Omega (U_3,L)\leq
  U_3=\Omega (N_1,L)$, we conclude $\norm(U_3)\leq N_3$; and hence,
  $S_4=\norm(U_3)$.
    
  Fourthly, $U_4=\Omega (N_2,U_3)$. Since
  \[U_3=(\PROD L)\langle\sigma\rangle (\DIAG U_3)\leq U_4\leq N_2
  =(\PROD U_2)\langle\sigma\rangle (\DIAG N_1)
  \]
  and $(\PROD U_2) \leq \Omega (N_2,U_2)\leq \Omega (N_2,U_3)$, we
  have $(\PROD U_2) \langle\sigma\rangle (\DIAG U_3)\leq U_4\leq
  N_2=(\PROD U_2) \langle\sigma\rangle (\DIAG N_1)$. As in the
  previous case, $U_4=(\PROD U_2) \langle\sigma\rangle (\DIAG V)$,
  with $V=\Omega(N_1,U_3)$.
    
  Finally, $N_4/N_3\cong (U_4/U_3)(S_4/S_3)$, and
  \[S_4/S_3=\norm(U_3)/N_2\cong((\PROD W)\langle\sigma\rangle (\DIAG N_2))
  /((\PROD U_2)\langle\sigma\rangle (\DIAG N_1)) \cong (W/U_2)(N_2/N_1).\]
\end{proof}

It was shown in~\cite{sidki:summernotes} that, up to conjugation, the
maximal $2$-subgroups $M$ of $\aut(\T)$ can be of two types: first,
subgroups of the form $M=M_1\times M_2$ where $M_1,M_2$ are two
nonconjugate maximal $2$-subgroups of $\aut(\T)$. The second type
consists of layered groups $M=(M\times M)\langle\sigma\rangle$.

\begin{cor}
  Suppose $L$ is a maximal $2$-subgroup of $\aut(\T)$ which is layered.
  Then $N_2(L)=N_1(L)$, i.e.\ $L$ has height at most $2$.
\end{cor}

We now proceed towards the proof of Theorem~\ref{thm:mainL}.
\begin{lem}\label{lem:Lrist}
  Let $L$ be a layered group. Then for every $i\in\N$ we have
  $\rist_{\aut^i(L)}(w)=L$ for all $w\in\T$ with $|w|\ge i$.

  The subgroup $\PROD^i L$ is characteristic in $\aut^i(L)$.
\end{lem}
\begin{proof}
  For $i=0$ the assertions are clear, while for $i=1$ they follow from
  Lemma~\ref{lem:norm}. We then prove both claims simultaneously by
  induction.
  
  Assume that $\rist_{\aut^i(L)}(w)=L$ for all $w\in\T$ with $|w|\ge
  i$, and $\PROD^i L\charin\aut^i(L)$. It follows from
  Theorem~\ref{thm:ln} that $\PROD^{i+1} L$ is characteristic in
  $\aut^i(L)$, and from Lemma~\ref{lem:norm} we have
  $\rist_{\aut^{i+1}(L)}(w)\le\aut(L)$ as soon as $|w|\ge i$. Therefore
  $\rist_{\aut^{i+1}(L)}(wx)\le\rist_{\aut(L)}(x)\le L$ for all $x\in\T$, and
  $\rist_{\aut^{i+1}(L)}(w)=L$ as soon as $|w|\ge i+1$.
  
  Then, again by Theorem~\ref{thm:ln},
  $\PROD^{i+1}L\charin\aut^{i+1}(L)$, because the product of all rigid
  stabilisers on a level of the tree is always characteristic.
\end{proof}

\begin{lem}\label{lem:wrdecomp}
  Let $L$ be a layered group, and set $N=\aut(L)$. Then $\aut^i(L)\le
  N\wr(\wr^{i-1} C_2)$ for all $i\in\N$.
\end{lem}
\begin{proof}
  By Theorem~\ref{thm:ln}, the terms $N_i$ of the automorphism tower
  of $L$ are subgroups of $W$. By Lemma~\ref{lem:Lrist}, the rigid
  stabilisers at level $i$ in $\T$ are all $L$. Given
  $\alpha\in\aut^{i+1}(L)$, we decompose it $i$ times as
  $\alpha=(\alpha_u,\dots,\alpha_v)\pi$ with $\alpha_w\in W$ for all
  $w\in\T$ of length $i$, and $\pi\in\wr^i C_2$.
  
  Since $\alpha$ normalizes $\aut^i(L)$, each $\alpha_u$ must
  normalize $\rist_{\aut^i(L)}(u)=L$, so $\alpha_u\in N$.
\end{proof}

\begin{proof}[Proof of Theorem~\ref{thm:mainL}]
  The first part follows directly from Lemma~\ref{lem:wrdecomp}.
  
  For the second part, Theorem~\ref{thm:normL} tells us how to compute
  $\norm^{i+1}(L)$ from $\norm^i(L)$. It remains to check that, given
  $\norm^i=(\PROD U_i)\langle\sigma\rangle(\DIAG S_i)$ with
  $U_i,S_i\in\mathcal C$, we have $U_{i+1}$ and $S_{i+1}\in\mathcal
  C$. This is clealy true for $U_{i+1}$, and for $S_{i+1}$ the
  computations of Theorem~\ref{thm:layered} show that it holds for
  $S_{i+1}$ if $i\le 3$.
  
  Consider then $i\ge 4$, and the group $U_i\in\mathcal C$. By
  Lemma~\ref{lem:Lrist}, we may write
  \[U_i=\left(\prod_{Z\in\{\PROD,\DIAG\}^i}Z\,G_Z\right)(\wr^i\langle\sigma\rangle),\]
  and we have $G_{\PROD^i}=L$.  The recursive computation of $\norm
  U_i$, following Theorem~\ref{thm:normL}, asks us to compute
  $\norm(G_{\PROD^i})=N$, and its intersection with some of the
  $\norm(G_Z)$; for this last step it suffices to compute
  $\norm_N(G_Z)$ and take intersections, all operations that leave
  elements of $\mathcal C$.

  The same argument applies to $\norm_{U_i}(S_i)=S_{i+1}$.
\end{proof}

We next construct within $W$ a family of layered subgroups of
arbitrary height $n\in\N$. Set
\[L_n=\langle\gamma_n(W),F\rangle,\]
where $\gamma_n(W)$ is the $n$th term in the lower central series of
$W$.
\begin{thm}
  $L_n$ is a layered group of height $n-1$; we have
  $\aut^m(L_n)=L_{n-m}$ for all $m\in\{0,\dots,n-1\}$.
\end{thm}
\begin{proof}
  We note that $L_n$ is layered by Lemma~\ref{lem:pre}.
  
  We recall the following facts from~\cite{bartholdi:lcs}: the
  quotient $\gamma_i(W)/\gamma_{i+1}(W)$ is isomorphic to the
  $\Z/2$-module $\{f:\N\to\Z/2\}$, and maps to
  $\gamma_{i+1}(W)/\gamma_{i+2}(W)$ with finite kernel and co-kernel
  by commutation with $\PROJ^k\sigma$, where $2^k||i$.
  
  It suffices to check $\norm L_{n+1}=L_n$ for all $n\in\N$. First,
  $[L_n,L_{n+1}]\le L_{n+1}$ because
  $[\gamma_n(W),F]\le\gamma_{n+1}(W)$ and
  $[\gamma_n(W),\gamma_{n+1}(W)]\le\gamma_{n+1}(W)$; therefore
  $\norm L_{n+1}\ge L_n$.
  
  Second, given $g\not\in L_n$, there exists $i<n$ with $g\in
  L_i\setminus L_{i+1}$; then by~\cite{bartholdi:lcs} there exists
  $f\in F$ with $[g,f]\not\in L_{i+2}$, so in particular $[g,f]\not\in
  L_{n+1}$ and $g\not\in\norm L_{n+1}$.
\end{proof}


\subsection{The dihedral, affine and cardioid group} We compute here the
automorphism and normalizer tower of three examples; the first two are
not layered, while the third one is.  We start by fixing some common
notation. Consider the adding machine $\tau=(1,\tau)\sigma$ acting on
the binary tree $\T=\{\6,\7\}^*$. It generates an infinite cyclic
group. Consider also the order-$2$ isometry $\delta$ of $\T$, defined
by $\delta=(\delta,\delta)\sigma$.
\begin{thm}
  The group $H=\langle\tau,\delta\rangle$ is infinite dihedral, and
  satisfies $\norm_W(H)=H$.
\end{thm}
\begin{proof}
  We have $\delta=u_{-1}\tau^{-1}$, and therefore
  $\tau^\delta=\tau^{-1}$ and $H$ is infinite dihedral.
  
  Then $\aut(H)$ is also dihedral, generated by $H$ and a new element
  $\upsilon$ satisfying $\upsilon^2=\tau$ and
  $\upsilon^\delta=\upsilon^{-1}$. We have $[\aut(H):H]=2$.
  
  However, there exists no tree isometry $g$ satisfying $g^2=\tau$;
  therefore $\norm(H)=H$ and the tower of $H$ strictly contains the
  normalizer tower.
\end{proof}
Note that $H$ is a finitely generated solvable group with $\norm(H)$
countable. In contrast, if $B$ is an abelian or finite $2$-subgroup of
$\aut\T$, then $\cent_{\aut\T}(B)$ is uncountable.

Consider next the topological closure of $\langle\tau\rangle$. It is
$T=\setsuch{\tau^x}{x\in\Z_2}\cong\Z_2$, the \emph{triangular group}.
Concretely, $\tau^x$ can be defined as
\[\tau^x=\begin{cases}
  (\tau^{x/2},\tau^{x/2}) & \text{ if }x\in 2\Z_2,\\
  (\tau^{(x-1)/2},\tau^{(x+1)/2})\sigma & \text{ if }x\in 2\Z_2+1.
\end{cases}
\]
Consider also the group
$U=\setsuch{u_x}{x\in\Z_2^\bullet}\cong\Z_2^\bullet$, the
\emph{unitary group}, where
\[u_x=(u_x,u_x\tau^{(x-1)/2}).\]
Set finally $A=TU$, the \emph{affine group}. The element $\tau^xu_y$
can be represented by the matrix
$(\begin{smallmatrix}1&x\\0&y\end{smallmatrix})$ over $\Z_2$.

\begin{thm}
  We have
  \[\norm_W(T)=A,\qquad\norm_W(A)=A.\]
\end{thm}
\begin{proof}
  We first note that $\cent_W(T)=T$. If $\alpha\in W$ fixes $\tau$,
  write $\alpha=(\alpha_\6,\alpha\7)\sigma^i$. Then
  $\tau^\alpha=\sigma^i(\alpha_\6^{-1}\alpha_\7,\alpha_\7^{-1}\tau\alpha_\6)\sigma^{1-i}=(1,\tau)\sigma$,
  so either $i=0,\alpha_\6=\alpha_\7\in\cent_W(\tau)$ or
  $i=1,\alpha_\6=\alpha_\7\tau^{-1}\in\cent_W(\tau)$; so in all cases,
  $\alpha\in T$.  Since $T$ has $\langle \tau\rangle$ as a dense
  subgroup, the claim follows.
  
  Now since $W$ acts continuously on itself by conjugation, any
  $a\in\norm_W(T)$ is determined by its value on $\tau$, and a simple
  computation shows that $\tau^{u_x}=\tau^x$, proving the first equality.
  
  Next, $T$ is a $\norm_W(A)$-invariant subgroup of $A$, since it can
  be defined as the set of fixed-point-free permutations of $A$;
  therefore any $a\in\norm_W(A)$ is again determined by its value on
  $\tau$, and the second equality follows.
\end{proof}
However, the automorphism tower of $T$ grows higher: even restricting
to continuous automorphisms, we have
$\aut^n(A)=\frac1{2^n}\Z_2\rtimes\Z_2^\bullet\cong A$, and
$\aut^\omega(A)=\Q_2\rtimes\Z_2^\bullet$, and
$\aut^{\omega+1}(A)=\aut^{\omega+2}(A)=\Q_2\rtimes\Q_2^\bullet$.

Our last example is the group $G=\langle F,\tau\rangle$, where $F$
denotes the ``finitary group'' generated by
$\setsuch{\6^n*\sigma}{n\in\N}$. We call $G$ the \emph{Cardioid
  group}; its name derives from the ``monodromy construction'' of
V.~Nekrashevych~\cite{bartholdi-g-n:fractal}, when applied to the
polynomial $f(z)=z^2+\epsilon$, with $\epsilon$ in the Mandelbrot
set's main cardioid. We briefly summarize its construction: let $M$ be
the Riemann surface $\C\setminus\overline{\{f^n(\epsilon)\}_{n\ge0}}$,
and let $M_0=M\setminus\{0\}$ be an open subsurface of $M$. Then $f$
defines a $2$-fold covering from $M_0$ to $M$. Let $t\in M$ be an
arbitrary point, and let $\T=\bigsqcup_{n\ge0}f^{-n}(t)\times\{n\}$ be
the disjoint union of all preimages of $t$. It has naturally the
structure of a binary tree, with edges connecting $(x,n+1)$ to
$(f(x),n)$. We let $\pi_1(M,t)$ act by monodromy (path lifting) on
$\T$, and denote by $\mathsf{IMG}(f)$ the quotient of $\pi_1(M,t)$ by
the kernel of the action.

It can be shown that $G=\mathsf{IMG}(f)$.  Indeed $\pi_1(M,t)$ is
free, and as generators of $\pi_1(M,t)$ one may choose a loop around
$\infty$, giving $\tau$, and loops encircling
$\overline{\{f^m(\epsilon)\}_{m\ge n}}$, giving $\6^n*\sigma$.
\begin{thm}\begin{enumerate}
  \item $G$ is a countable layered group, and satisfies
    $G/G'=\Z\oplus(\Z/2)^\infty$.
  \item $\norm(G)=\langle F,T\rangle D$, where $D=\langle
    \delta\rangle\cong\Z/2$.
  \item $\norm^2(G)=\norm(G)$.
  \end{enumerate}
\end{thm}
\begin{proof}
  Clearly $G$ is layered, because $\tau\sigma=(1,\tau)\in G$.  The
  abelianization of $W$ is $(\Z/2)^\N$, and the image of $G$ in $W/W'$
  is generated by the finite-support and constant functions.
  Furthermore, the image of $\tau$ in $G/G'$ has infinite order, since
  $\tau^{2^n}\in G'$ implies $\tau\in G'$ contradicting the above.
  
  For the second assertion, we first check that $T$ normalizes $F$;
  indeed
  \[\sigma^{\tau^x}=\begin{cases}(\tau^{-x/2},\tau^{-x/2})\sigma(\tau^{x/2},\tau^{x/2}) & \text{ if }x\in2\Z_2\\
    \sigma(\tau^{(1-x)/2},\tau^{(-1-x)/2})\sigma(\tau^{(x-1)/2},\tau^{(x+1)/2})\sigma & \text{ if }x\in2\Z_2+1\end{cases}=\sigma,\]
  and then by induction
  \[(\6^n*\sigma)^{\tau^x}=\begin{cases}((\6^{n-1}*\sigma)^{\tau^{x/2}},1) &
  \text{ if }x\in2\Z_2\\
    (1,(\6^{n-1}*\sigma)^{\tau^{(x-1)/2}}) & \text{ if }x\in2\Z_2+1\end{cases}\in\sigma^F.
  \]
  Since $T$ clearly centralizes $\tau$, it follows that $T$ normalizes
  $G$.
  
  Next, we compute $\norm(G)\cap A=TD$; indeed
  $\sigma^{u_x}=(\tau^{(x-1)/2},(\tau^{(x-1)/2})^{u_x})\sigma$ belongs to
  $G$, and therefore $(x-1)/2\in\Z$, so $x\in2\Z+1$ is invertible, and
  $x=\pm1$.

  Finally, consider an arbitrary $a\in\norm(G)$. Since $\tau^a\in G$, we
  have $\tau^a=(\tau^{i_{\6\dots\6}},\dots,\tau^{i_{\7\dots\7}})\pi$ for some
  permutation $\pi$ of the first $n$ levels of $\T$. Furthermore
  $\pi$ is $2^n$-cycle, and therefore $(\tau^{2^n})^a=(\tau^j,\dots,\tau^j)$
  for some odd $j\in\Z$.
  
  On the other hand, since $G$ is layered, we may decompose
  $a=(b,\dots,b)(1,g_{\6\dots\6\7},\dots,g_{\7\dots\7})\pi$ for some
  $g_v\in G$ and $b\in\norm(G)$. Then $\tau^b=\tau^j$, and therefore
  $b\in TD$, so $a\in FTD$.
  
  The last assertion follows from Theorem~\ref{thm:layered}.  Indeed
  in that theorem's notation $U_2=\Omega(\norm(G),G)=G$, since
  $[\tau^xd,\tau^y]=\tau^{2y}\not\in G$ as soon as $y\not\in\Z$, whence
  $U_2\le FT$; and $(\tau^x)^2\not\in G$ for $x\not\in\Z$, whence
  $U_2=F\langle \tau\rangle$. Therefore
  $\norm^2(G)=\norm(G)(\PROD U_2)=\norm(G)$.
\end{proof}

\section{The three groups $\GS,\BG,\FG$ acting on the ternary tree}
\label{sec:gs}
We consider in this section groups acting on the ternary tree $\T=X^*$
with $X=\{\6,\7,\8\}$, that are generated by two elements $x$ and
$\gamma$. The $3$-cycle $x=(\6,\7,\8)$ induces a rigid permutation of
the subtrees $\6\T,\7\T,\8\T$ of $\T$ as in~\eqref{eq:rigid}. Choose
now $i,j\in\{0,1,2\}$, and consider the element
$\gamma=(\gamma,x^i,x^j)$ acting recursively as $\gamma,x^i,x^j$ on
the respective subtrees $\6\T,\7\T,\8\T$. Up to conjugation by a tree
isometry, it is sufficient to consider such groups with $i=1$, and
there are therefore at most three non-isomorphic such groups. These
groups are known to be non-isomorphic~\cite{bartholdi-g:parabolic}. We
concentrate here on the first example, the group $\GS$ for which
$\gamma=(\gamma,x,x^{-1})$. It was first studied
in~\cite{gupta-s:burnside}.

Recall that a group is \emph{just-infinite} if it is infinite, but all
its proper quotients are finite.  The main properties of $\GS$ are:
\begin{thm}
  $\GS$ is a just-infinite $3$-group, and it is regular branch: we
  have $\GS'>\PROD\GS'$ with $\GS/\GS'\cong(\Z/3)\times(\Z/3)$ and
  $\GS'/(\GS'\times \GS'\times \GS')\cong(\Z/3)\times(\Z/3)$.
\end{thm}

\subsection{The original strategy}
This section describes the automorphism tower of the group $\GS$, its
normalizer in the pro-$3$-Sylow of $W$, and sketches similar results
for two related groups acting on the ternary tree. We start by
reviewing the strategy followed by the second author
in~\cite{sidki:subgroups} to compute the automorphism group of $\GS$:

\begin{enumerate}
\item Define the ascending series of subgroups $\GS^0=\GS$ and
  $\GS^{n+1}=\GS^n\wr\langle x\rangle$; all these groups act
  level-transitively on $\T$.  Set $\GS_*=\bigcup_{n\in\N}\GS^n$.
\item Observe that $\rist_\GS(\6)\wr\langle x\rangle\leq
  \GS\leq\GS\wr\langle x\rangle=\GS^1$. Since
  $\GS^1/(\rist_\GS(\6)\wr\langle x\rangle)$ is a finite $3$-group, it
  follows that $\norm_{\GS^1}(\GS)\ge\GS$; check that
  $\norm_{\GS^1}(\GS) =\GS\langle(x,x,x)\rangle$.  This procedure is
  repeated for the pair $(\rist_\GS\wr^n\langle x\rangle,\GS^n)$ for all
  $n\in\N$, and yields $\norm_{\GS_*}(\GS)$.
\item Determine $3'$-automorphisms of $\GS$: produce a Klein $4$-group
  acting simply on the generating set
  $\{x,x^{-1},\gamma,\gamma^{-1}\}$ of $\GS$.
\item The problem is reduced to the description of the
  \textsf{IA}-automorphisms\footnote{i.e.\ those automorphisms which
    induce the identity on $\GS/\GS'$} of $\GS$.  Express the
  \textsf{IA}-automorphisms as (potentially infinite) words in
  $\norm_{\GS_*}(\GS)$.
\item Eliminate the possibility of infinite words (i.e.\ elements of
  $\widehat{\norm_{\GS_*}(\GS)}$) by using induction on the depth
  function of $\GS$, defined on $g$ as the minimal level $n$ at which
  all states $g@v$ belong to $\{1,x^{\pm1},\gamma^{\pm1}\}$ for all
  $v\in X^n$.
\end{enumerate}

\begin{thm}[\cite{sidki:subgroups}]
  Let $t$ be the transposition $(\7,\8)$. Set
  $U=\langle\DIAG^n(x):n\ge1\rangle$, a countably infinite
  elementary abelian $3$-group, and $V=\{1,\tau_1,\tau_2,\tau_3\}$, a
  Klein group, where
  \[\tau_1=\DIAG(\tau_1)t,\quad\tau_2=\DIAG(\tau_3),\quad\tau_3=\DIAG(\tau_2)t.\]
  Then we have a split extension
  \[\aut(\GS)=\GS:U:V.\]
  Furthermore $V$ normalizes $U$, and acts on it via
  \begin{xalignat*}{3}
    \DIAG^n(x)^{\tau_1}&=\DIAG^n(x)^{-1},&
    \DIAG^n(x)^{\tau_2}&=\DIAG^n(x)^{(-1)^n},&
    \DIAG^n(x)^{\tau_3}&=\DIAG^n(x)^{-(-1)^n},\\
    \gamma^{\tau_1}&=\gamma,&
    \gamma^{\tau_2}&=\gamma^{-1},&
    \gamma^{\tau_3}&=\gamma^{-1}.
  \end{xalignat*}
\end{thm}

\begin{thm}\label{thm:gstower}
  Set $T=\langle \DERI\DIAG^{2n+1}(x):n\ge0\rangle$, a countably
  infinite elementary abelian $3$-group. Then
  \[\aut^2(\GS)=\aut(\GS):T.\]
  Furthermore $T$ centralizes $U$ and $V$.

  \[\aut^3(\GS)=\aut^2(\GS).\]
\end{thm}
\begin{rem}  
  It follows from the proof that $\GS'$ is absolutely characteristic
  in $\GS$, and $\aut^2(\GS)$ is a subgroup of the commensurator of
  $\GS$.
\end{rem}

To simplify the notations, we write $A=\aut(\GS)$ and $B=\aut^2(\GS)$. For
$v\in X^n$, we also write
\[\stab_\GS^v(n)=\setsuch{g\in \stab_\GS(n)}{g@v=1}.\]

\begin{lem}[\cite{sidki:subgroups}]
  $\GS'$ is a subdirect product of $\GS\times \GS\times \GS$, and factors as
  \[
  \GS'=\left(\GS'\times \GS'\times \GS'\right)\langle[\gamma,x],\DIAG\gamma\rangle.
  \]
  Furthermore, $\cent_W(\GS')$ is trivial.
\end{lem}

\begin{lem}
  $A'=\GS:U$, and $A''=\GS'$. Therefore $\GS'$ is characteristic in $A$.
\end{lem}
\begin{proof}
  The first assertion follows immediately from the fact that the every
  element of the generating set
  $\left\{\gamma,\DIAG^n(x):n\ge0\right\}$ of $\GS U$ is inverted by
  $\tau_1$ or $\tau_2$. The second assertion follows from
  $A''=\GS'[\GS,U]$, and $[x,\DIAG^n(x)]=1$,
  $[\gamma,\DIAG^n(x)]=([\gamma,\DIAG^{n-1}(x)],1,1)\in \GS'$.
\end{proof}

\begin{lem}
  $B\le(A\times A\times A)\langle x,\tau_1\rangle$.
\end{lem}
\begin{proof}
  First, note that $A$ is saturated (the subgroups $\langle
  v*\GS':|v|=n\rangle$ may serve as $H_n$) and weakly branch (since it
  contains $\GS$), and apply Theorem~\ref{thm:ln} to obtain
  $B\cong\norm_W(A)$; these two groups will be identified. As $\GS'$
  is $B$-invariant and $\cent_W(\GS^{\prime})=\{e\}$, we conclude that
  $B\le\aut(\GS')$.
  
  Given $\phi\in B$, it can be written as
  $\phi=(\phi_\6,\phi_\7,\phi_\8)\tau_1^ix^j$ for some $i\in\{0,1\}$
  and $j\in\{0,1,2\}$.  We conclude $(\phi_\6,\phi_\7,\phi_\8)\in B$.
  Now since $\GS'$ is a subdirect product of $\GS\times \GS\times \GS$, we
  obtain $\phi_x\in\norm_W(\GS)$ for $x\in X$. The conclusion follows.
\end{proof}

\begin{lem}[\cite{sidki:subgroups}]
  Take $g\in \GS$. Then $g\in\stab_\GS(1)$ if and only if
  \[g=\left(c_\6\gamma^{n_1}x^{n_3-n_2},c_\7\gamma^{n_2}x^{n_1-n_3},c_\8\gamma^{n_3}x^{n_2-n_1}\right)\]
  for some $c_x\in \GS'$ and $n_1,n_2,n_3\in\{0,1,2\}$. We have
  \begin{align*}
    \stab_\GS(1)&=\left(\GS'\times \GS'\times \GS'\right)
    \langle\DIAG\gamma,(1,\gamma x^{-1},\gamma
    x),(x,x^{-1},\gamma)\rangle,\\
    \text{and }\stab_\GS^\6(1)&=\left(1\times \GS'\times \GS'\right)\langle(1,\gamma x^{-1},\gamma x).
  \end{align*}
\end{lem}

\begin{lem}
  \begin{align*}
    \stab_A(1)&=\stab_\GS(1)U\langle\tau_2\rangle,\\
    \text{and }\stab_A^\6(1)&=\stab_\GS^\6(1)\langle(1,\gamma^{-1},\gamma)\rangle,
  \end{align*}
  which is a subdirect product of $1\times \GS\times \GS$.
\end{lem}
\begin{proof}
  The first statement is clear.  Then take $\phi\in\stab_A^\6(1)$; we
  can write
  \[\phi=\left(c_\6\gamma^{n_1}x^{n_3-n_2}uv,c_\7\gamma^{n_2}x^{n_1-n_3}uv,c_\8\gamma^{n_3}x^{n_2-n_1}uv\right)\]
  for some $u\in U@\6=U\langle x\rangle$ and $v\in V$. Since
  $c_\6\gamma^{n_1}x^{n_3-n_2}uv=1$, we have $n_1=0,c_\6=1,v=1$, and
  $u\in\langle x\rangle$ with $x^{n_3-n_2}u=1$. Modulo
  $\stab_\GS^\6(1)$, it suffices to consider $u=x$ and $n_2=2,n_3=1$,
  leading to $(1,\gamma^{-1},\gamma)$.
  
  The last assertion follows, since the second and third coordinates
  of $\langle(1,\gamma x^{-1},\gamma x),(1,\gamma^{-1},\gamma)\rangle$
  each generate $\GS$.
\end{proof}

\begin{lem}\label{lem:gsrw}
  Let $H$ be a group such that $A\le H\le(A\times A\times A)\langle
  x\rangle V$, and take $\phi\in H$. Then $\phi$ can be reduced modulo
  $A$ to $(v_\6,u_\7v_\7,x^i\gamma^ju_\8v_\8)$.
\end{lem}
\begin{proof}
  Direct.
\end{proof}

\begin{lem}
  Consider $\phi=\DERI\DIAG^{2n+1}(x)$, for $n\ge0$. Then $\phi\in
  B$ and $[UV,\phi]=1$.
\end{lem}
\begin{proof}
  To show that $\phi\in B$ it suffices to show $[g,\phi]\in A$ for all
  generators $g$ of $A$. First, $[x,\phi]=\DIAG^{2n+2}(x)$ and
  $[\gamma,\phi]=[\DIAG^m(x),\phi]=1$. Then we have
  \begin{align*}
    [\tau_1,\phi]&=\DERI\DIAG^{2n+1}(x)^{-\tau_1}\phi=\DERI(\DIAG^{2n+1}(x)^{\tau_1})\phi=\DERI\DIAG^{2n+1}(x)^{-1}\phi=1,\\
    \text{and }[\tau_2,\phi]&=\DERI\DIAG^{2n+1}(x)^{-\tau_2}\phi=\DERI(\DIAG^{2n+1}(x)^{-\tau_3})\phi=\DERI\DIAG^{2n+1}(x)^{-1}\phi=1,
  \end{align*}
  since $\tau_1$ inverts $\DIAG^{2n+1}(x)$ and $\tau_2$ fixes
  $\DIAG^{2n+1}(x)$.
\end{proof}

\begin{proof}[Proof of Theorem~\ref{thm:gstower}, first part]
  Take $\phi$ normalizing $A$, and write it, using
  Lemma~\ref{lem:gsrw}, as
  $\phi=(v_\6,u_\7v_\7,x^i\gamma^{j}u_\8v_\8)$. We compute the
  commutator of $\phi$ with generators of $A$ in turn:
  \begin{description}
  \item[interaction with $\tau_2$]
    \begin{align*}
      \phi^{\tau_2} &= (v_\6,u_\7v_\7,x^i\gamma^{j}u_\8v_\8)^{\DIAG\tau_3}=(v_\6,u_\7^{\tau_3}v_\7,x^{-i}\gamma^{-j}u_\8^{\tau_3}v_\8),\\
      [\tau_2,\phi^{-1}] &=
      (1,[\tau_3,u_\7^{-1}],x^{-i}\gamma^{-j}[\tau_3,u_\8^{-1}]\gamma^{-j}x^{-i})\in\stab_A^\6(1).
    \end{align*}
    Therefore $[\tau_3,u_\7^{-1}]=[\tau_3,u_\8^{-1}]=e$ and $i=j=0$;
    hence $\phi=(v_\6,u_\7v_\7,u _\8v_\8)$, where $u_\7,u_\8\in
    \cent(\tau_3)$; that is, are products of $\DIAG^{2n+1}(x)$. Now
    using $T$ we may reduce $\phi$ further to the form
    $\phi=(v_\6,v_\7,uv_\8)$, for some $u\in U$ a product of
    $\DIAG^{2n+1}(x)$.
  \item[interaction with $\tau_1$]
    \begin{align*}
      \phi^{\tau_1} &= (v_\6,v_\7,uv_\8)^{\DIAG(\tau_1)t}=(v_\6,uv_\8,v_\7)^{\DIAG\tau_1}=(v_\6,u^{-1}v_\8,v_\7),\\
      [\tau_1,\phi^{-1}] &= (v_\6,u^{-1}v_\8,v_\7)(v_\6,v_\7,(u
      v_\8)^{-1})=(1,u^{-1}v_\8v_\7,v_\7v_\8u^{-1})\in\stab_A^\6(1). 
    \end{align*}
    Therefore $u=1$, $v_\7=v_\8$, and $\phi$ is reduced to the form
    $\phi=(v_\6,v_\7,v_\7)$.
  \item[interaction with $x$]
    \[
      [x,\phi^{-1}]^x=(1,v_\7v_\6,v_\6v_\7)\in\stab_A^\6(1).
    \]
    Therefore $v_\6=v_\7$, and $\phi$ is reduced to $\phi=(v,v,v)$.
    If $v=\tau_3$ we have $\phi=\tau_2\in A$, so it suffices to
    consider $\phi=\DIAG\tau_1$.
  \item[interaction with $\gamma$]
    \[
      [\DIAG\tau_1,\gamma^{-1}]=(\gamma,x^{-1},x)(\gamma,x,x^{-1})^{-1}=(1,x,x),
    \]
    but this last element is not in $\stab_A^\6(1)$; therefore, this
    last case is impossible, and these successive steps have reduced
    $\phi$ to $1$.
  \end{description}
\end{proof}

\begin{lem}
  $B'=\GS:U=A'$, and $B''=\GS'=A''$.

  \noindent $\stab_B^\6(1)=\stab_A^\6(1):T$.

  \noindent $B$ is a saturated subgroup of $W$, and $\aut B\le(A\times
  A\times A)\langle x,\tau_1\rangle$.
\end{lem}
\begin{proof}
  Direct.
\end{proof}

\begin{proof}[Proof of Theorem~\ref{thm:gstower}, second part]
  Let $\phi$ normalize $B$. As before, we may assume
  $\phi=(v_\6,u_\7v_\7,x^i\gamma^{j}u_\8v_\8)$.
  Now, 
  \[
  [\tau_2,\phi^{-1}]=(1,[\tau_3,u_\7^{-1}],x^{-i}\gamma^{-j}[\tau_3,u_\8^{-1}]\gamma^{-j}x^{-i})\in\stab_B^\6(1).
  \]
  Therefore $i=j=0$, and $(1,[\tau_3,u_\7^{-1}],[\tau_3,u_\8^{-1}])\in
  T$. On one hand, $T$ is centralized by $\tau_3$; but on the other
  hand, $[\tau_3,u_\7^{-1}]$ and $[\tau_3,u_\8^{-1}]$ are inverted by
  $\tau_3$, since $\tau_3$ is an involution. Since $T$ has no
  $2$-torsion, it follows that
  $[\tau_3,u_\7^{-1}]=[\tau_3,u_\8^{-1}]=1$. We then have
  $\phi=(v_\6,u_\7v_\7,u_\8v_\8)$ where $u_\7,u_\8$ are centralized by
  $\tau_3$. Using $(1,u_\7,u_\7^{-1})\in T$, we may reduce to the case
  $\phi=(v_\6,v_\7,u v_\8)$.

  We next compute
  \begin{align*}
    \phi^{\tau_1} &= (v_\6,(u v_\8)^{\tau_1},v_\7)=(v_\6,u^{-1}v_\8,v_\7),\\
    [\tau_1,\phi^{-1}] &=
    (1,u^{-1}v_\8v_\7,v_\7v_\8u^{-1})\in\stab_B^\6(1).
  \end{align*}
  Therefore $v_\7=v_\8$ and $[\tau_1,\phi^{-1}]=(1,u^{-1},u^{-1})\in
  UT$, hence $u=1$ and the proof is finished.
\end{proof}

\subsection{The normalizer tower}\label{subs:GStower}
The reason the automorphism tower of $\GS$ stops after $2$ steps can be
attributed to the involutions in $V$; just in the same way as the
lower central series of $\aut(\GS)$ stabilizes at $\GS:U$ and that of
$\aut^2(\GS)$ stabilizes at $\GS:UT$, the normalizer tower may not be
extended beyond $V$.

However, since $\GS$ is a subgroup of the pro-$3$-Sylow
$C=\widehat{\langle x\rangle_*}$ of $W$, we may consider its
normalizer tower within $C$.  For this purpose, consider the set
$\Omega=\setsuch{Z\,x}{Z\in\{\DIAG,\DERI,\PROJ\}^*}$, with the notation
of Definition~\ref{defn:xaction}. The set $\Omega$ will serve as a
basis for the quotients of successive terms of the normalizer tower of
$\GS$. Define a rank function on $\Omega$ by
\begin{xalignat*}{2}
  \rank x&=1,&\rank(\PROJ x)&=3\rank x,\\
  \rank(\DERI x)&=3\rank x-1,&\rank(\DIAG x)&=3\rank x-2,
\end{xalignat*}
and define $A_n=\langle x,\gamma,Z\,x\text{ such that }\rank(Z\,x)\le n\rangle$.
The elements of $\Omega$ of low rank are as given below:
\[\begin{array}{r|ccccccccc}
  \text{rank} & \multispan{2}\text{\qquad elements}\\ \hline
  1 & x & \DIAG x & \DIAG\DIAG x & \DIAG\DIAG\DIAG x & \DIAG\DIAG\DIAG\DIAG x & \cdots & \DIAG^nx\\
  2 &   & \DERI x & \DERI\DIAG x & \DERI\DIAG\DIAG x & \DERI\DIAG\DIAG\DIAG x & \cdots & \DERI\DIAG^nx\\
  3 &   & \PROJ x & \PROJ\DIAG x & \PROJ\DIAG\DIAG x & \PROJ\DIAG\DIAG\DIAG x & \cdots & \PROJ\DIAG^nx\\
  4 &   &         & \DIAG\DERI x & \DIAG\DERI\DIAG x & \DIAG\DERI\DIAG\DIAG x & \cdots & \DIAG\DERI\DIAG^nx\\
  5 &   &         & \DERI\DERI x & \DERI\DERI\DIAG x & \DERI\DERI\DIAG\DIAG x & \cdots & \DERI\DERI\DIAG^nx\\
  6 &   &         & \PROJ\DERI x & \PROJ\DERI\DIAG x & \PROJ\DERI\DIAG\DIAG x & \cdots & \PROJ\DERI\DIAG^nx\\
  7 &   &         & \DIAG\PROJ x & \DIAG\PROJ\DIAG x & \DIAG\PROJ\DIAG\DIAG x & \cdots & \DIAG\PROJ\DIAG^nx\\
  8 &   &         & \DERI\PROJ x & \DERI\PROJ\DIAG x & \DERI\PROJ\DIAG\DIAG x & \cdots & \DERI\PROJ\DIAG^nx\\
  9 &   &         & \PROJ\PROJ x & \PROJ\PROJ\DIAG x & \PROJ\PROJ\DIAG\DIAG x & \cdots & \PROJ\PROJ\DIAG^nx\\
 10 &   &         &              & \DIAG\DIAG\DERI x & \DIAG\DIAG\DERI\DIAG x & \cdots & \DIAG\DIAG\DERI\DIAG^nx
\end{array}\]
The next rows follow the same pattern: elements of rank $n$ are
determined by writing backwards $n-1$ in base $3$, using the symbols
$\DIAG=0$, $\DERI=1$ and $\PROJ=2$.

\begin{thm}
  The following facts hold for all $n\ge1$:
  \begin{enumerate}
  \item $\norm_C^n(\GS)\ge A_n$;
  \item each $\norm_C^{n+1}(\GS)/\norm_C^n(\GS)$ is an infinite countable
    elementary abelian $3$-group;
  \item $\norm_C^\alpha(\GS)=\langle\gamma,\Omega\rangle=\GS_*$ for
    $\alpha\ge\omega$.
  \end{enumerate}
\end{thm}
\begin{proof}
  Extend the rank function to $\Omega\cup\{\gamma\}$, by
  $\rank\gamma=0$.  Inductive computations then show that
  $\rank[u,v]\le\max\{\rank u,\rank v\}-1$ for all
  $u\in\Omega\cup\{\gamma\}$ and $v\in\Omega$ of rank at least $1$.
  For example, for any $Y,Z\in\{\DIAG,\DERI,\PROJ\}^*$:
  \begin{align*}
    \rank[\DERI Y x,\DIAG Z x]&=\rank\DERI[][Y x,Z x]\\
    &\le3\rank[Y x,Z x]-1\\
    &\le3(\max\{\rank(Y x),\rank(Z x)\}-1)-1\\
    &\le\rank\max\{\rank(\DERI Y x),\rank(\DERI Y x)\}-1,\\
\text{and }\rank[x,\DERI Z x]&=\rank\DIAG Z x=\rank(\DERI Z x)-1.
  \end{align*}
  We therefore have $A_n\triangleleft A_{n+1}$ with abelian quotients,
  so the tower $A_n$ grows more slowly than the normalizer tower.

  Finally, the tower of $\GS$ stops at level $\omega$ by
  Theorem~\ref{thm:stops@omega}.
%
\end{proof}

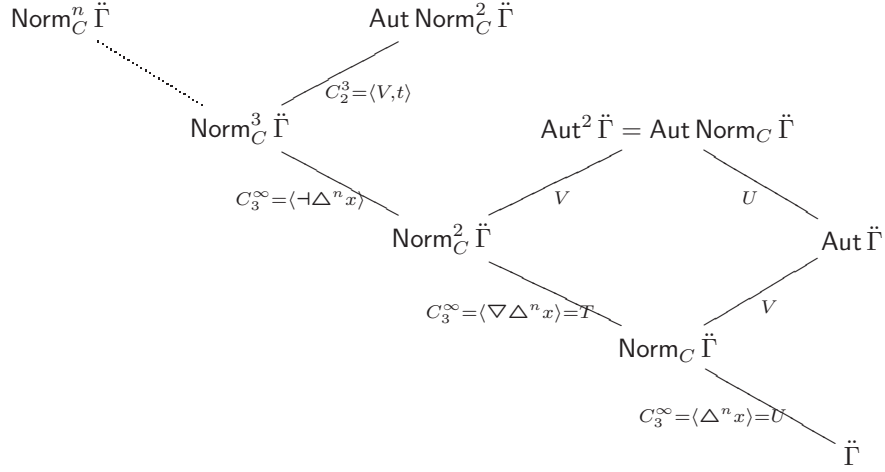
\[\xymatrix{{\norm^n_C\GS}\ar@{.}[dr] & &
  {\aut\norm^2_C\GS}\ar@{-}[dl]^{C_2^3=\langle V,t\rangle}\\
  & {\norm^3_C\GS}\ar@{-}[dr]_{C_3^\infty=\langle\PROJ\DIAG^nx\rangle}
  & & {\kern-2em\aut^2\GS=\aut\norm_C\GS\kern-2em}\ar@{-}[dl]^{V}\ar@{-}[dr]_{U}\\
  & & {\norm^2_C\GS}\ar@{-}[dr]_{C_3^\infty=\langle\DERI[\scriptstyle]\DIAG^nx\rangle=T}
  & & {\aut\GS}\ar@{-}[dl]^{V}\\
  & & & {\norm_C\GS}\ar@{-}[dr]_{C_3^\infty=\langle\DIAG^nx\rangle=U}\\
  & & & & {\GS}
}\]

\subsection{The Fabrykowski-Gupta group}
This group was introduced in~\cite{fabrykowski-g:growth1} as an
alternate construction of a group of intermediate growth. Like $\GS$,
it acts on the ternary tree $\{\6,\7,\8\}^*$. Let $x=(\6,\7,\8)$ be
the $3$-cycle acting on $\T$ according to~\eqref{eq:rigid}, and define
recursively $\gamma=(\gamma,x,1)$. The Fabrykowski-Gupta group is
$\FG=\langle x,\gamma\rangle$. It is a regular branch group; we have
$\FG/\FG'\cong(\Z/3)\times(\Z/3)$, and $\FG'/(\FG'\times \FG'\times
\FG')\cong(\Z/3)\times(\Z/3)$.
  
The generators $x,\gamma$ of $\FG$ have order $3$, and $\langle
x\gamma,\gamma x\rangle$ is a torsion-free normal subgroup of $\FG$ of
index $3$.
  
\begin{thm}
  Set $U=\langle\DIAG^n(x):n\ge1\rangle$, a countably infinite
  elementary abelian $3$-group. Then we have a split extension
  \[\aut(\FG)=\FG:U.\]
\end{thm}
\begin{proof}
  The proof follows the same scheme as that of
  Theorem~\ref{thm:gstower} and is omitted.
\end{proof}

\subsection{The group $\BG$}
This group also acts on the ternary tree, and is generated by
$x=(\6,\7,\8)$ and $\gamma=(\gamma,x,x)$.  The group $\BG$ was first
studied in~\cites{bartholdi-g:spectrum,bartholdi-g:parabolic}, where
some of its elementary properties were proven.

Note that the group $\BG$ is not branch, but is a regular weakly
branch group. Consider the subgroup $K=\langle x^{-1}\gamma,\gamma
x^{-1}\rangle$; then $K$ is a torsion-free normal subgroup of index
$3$. We have $\BG/\BG'\cong(\Z/3)\times(\Z/3)$, and
$\BG'/K'\cong\Z^2$, and $K'/(K'\times K'\times K')\cong\Z^2$.
  
\begin{thm}
  Let $t$ be the transposition $(\7,\8)$ acting on $\T$ according
  to~\eqref{eq:rigid}. Let
  $U=\langle\DIAG^n(x\gamma^{-1}):n\ge1\rangle$ be a countably
  infinite elementary abelian $3$-group, and let
  $V=\{1,\tau_1,\tau_2,\tau_3\}$ be the Klein group defined by
  \[\tau_1=\DIAG(\tau_1)t,\quad\tau_2=\DIAG(\tau_1),\quad\tau_3=t.\]
  Then we have a non-split extension
  \[\aut(\BG)=(\BG\cdot U):V,\]
  where $V$ acts on $U$, for $n>0$, by
  \begin{xalignat*}{3}
    x^{\tau_1}&=x^{-1},& x^{\tau_2}&=x,& x^{\tau_3}&=x^{-1},\\
    \DIAG^n(x)^{\tau_1}&=\DIAG^n(x)^{-1},&
    \DIAG^n(x)^{\tau_2}&=\DIAG^n(x)^{-1},&
    \DIAG^n(x)^{\tau_3}&=\DIAG^n(x),\\
    \gamma^{\tau_1}&=\gamma^{-1},&
    \gamma^{\tau_2}&=\gamma^{-1},&
    \gamma^{\tau_3}&=\gamma.
  \end{xalignat*}

  Let $T=\langle\DERI\DIAG^{2n+1}(x):n\ge0\rangle$ be a countably
  infinite elementary abelian $3$-group. Then
  \[\aut^3(\BG)=\aut^2(\BG)=\aut(\BG):T.\]
  Furthermore $T$ centralizes $U$ and $V$.
\end{thm}
\begin{proof}
  The proof follows the same scheme as that of
  Theorem~\ref{thm:gstower} and is omitted.
\end{proof}

\section{The First Grigorchuk Group}
\label{sec:grig}
This section describes automorphism tower of the first Grigorchuk
group. Consider the regular binary tree $X^*$ with $X=\{\6,\7\}$,
denote by $\sigma$ the non-trivial permutation $(\6,\7)$ of $X$ and
define isometries $a,b,c,d$ of $\T$ by
\begin{align*}
  (x_1\dots x_m)^a &= x_1^\sigma x_2\dots x_m,\\
  b &= (a,c) = \prod_{i\not\equiv 0\mod 3}\7^i\6*a,\\
  c &= (a,d) = \prod_{i\not\equiv 2\mod 3}\7^i\6*a,\\
  d &= (1,b) = \prod_{i\not\equiv 1\mod 3}\7^i\6*a.
\end{align*}
The Grigorchuk group~\cites{grigorchuk:burnside,grigorchuk:growth} is
$\Gg=\langle S\rangle$, with $S=\{a,b,c,d\}$.  It is a regular branch
group. Set $B=\rist_\Gg(\6)@\6=\langle b\rangle^G$ and
$K=\rist_\Gg(\6^2)@\6^2=\langle [a,b]\rangle^\Gg$.  Then $B$ has index
$8$ in $\Gg$, and $K$ has index $16$; and $K/(K\times K)\cong\Z/4$.

The automorphism group of $\Gg$ was computed
in~\cite{grigorchuk-s:auto}, where it was shown, in a manner similar
to the one employed in~\cite{sidki:subgroups}, that $\out(\Gg)$ is an
infinite elementary $2$-group, generated by $\DIAG^n\PROJ adad$ for
all $n\ge1$, in the notation of Definition~\ref{defn:xaction}. We will
give below a complete description of the tower, without relying on the
above result. In order to illustrate our method, we first examine more
closely the first terms of the automorphism tower.

Let us define inductively the normalizer tower by $N_0=\Gg$ and
$N_{n+1}=\norm_W(N_n)$. Then $N_1=\langle
a,b,c,d,\DIAG^n\PROJ adad\rangle$. There are clearly some elements of $W$
outside $N_1$ that normalize $N_1$; for instance, $\PROJ ad$. In fact, a
little computation shows that the elements $\PROJ\DIAG^n\PROJ adad$ for all
$n\in\N$ also normalize $N_1$: indeed for all $m,n\ge0$:
\begin{align*}
  [\PROJ\DIAG^n\PROJ adad,\DIAG^m\PROJ adad]&=\begin{cases}
    \PROJ\DIAG^{m-1}\PROJ[\DIAG^{n-m-1}\PROJ adad,\DIAG b]=1&\text{ if }m<n+1;\\
    1&\text{ if }m=n+1,\\
    \PROJ\DIAG^n\PROJ[\DIAG b,\DIAG^{m-n-2}\PROJ adad]=1&\text{ if }m>n+1;
  \end{cases}\\
  [\PROJ\DIAG^n\PROJ adad,a]&=\DIAG^{n+1}\PROJ adad\in N_1;\\
  [\PROJ\DIAG^n\PROJ adad,b]&=[\PROJ\DIAG^n\PROJ adad,c]=[\PROJ\DIAG^n
  adad,d]=1.
\end{align*}
Similarly, we have $\PROJ ad$ in $N_2$.

Proceeding, we see $\PROJ a$ in $N_3$, and therefore $N_3$ contains
$\PROD\Gg$. We may also check that $N_3$ contains some other
elements, like $\DIAG\PROJ ad$ and $\DIAG\PROJ\PROJ adad$. Based on the
behaviour prescribed by Theorem~\ref{thm:layered} for layered groups,
we may ``guess'' that $N_3$ is $(\PROD\Gg)(\DIAG N_1)\langle
a\rangle$.

Actually, the tower of automorphisms of $\Gg$ follows a pattern that
is best described by powers of $2$. We describe it in the same
language as the normalizer tower of $\GS$ given in
Subsection~\ref{subs:GStower}.

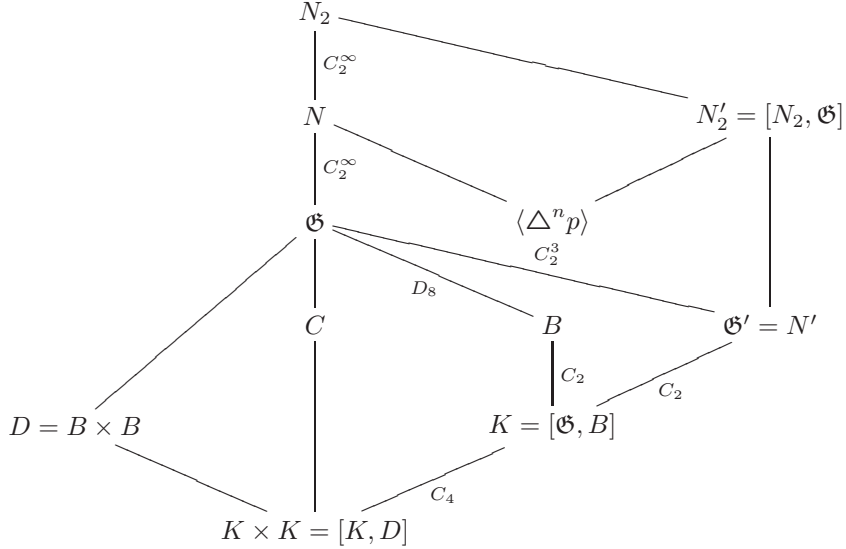
\[\xymatrix{ & {N_2}\ar@{-}[d]^{C_2^\infty}\ar@{-}[drr]\\
  & {N}\ar@{-}[d]^{C_2^\infty}\ar@{-}[dr] & & {N_2'=[N_2,\Gg]}\ar@{-}[dl]\ar@{-}[dd]\\
  & {\Gg}\ar@{-}[ddl]\ar@{-}[d]\ar@{-}[dr]_{D_8}\ar@{-}[drr]^{C_2^3} & {\langle\DIAG^np\rangle}\\
  & {C}\ar@{-}[dd] & {B}\ar@{-}[d]^{C_2} & {\Gg'=N'}\ar@{-}[dl]^{C_2}\\
  {D=B\times B}\ar@{-}[dr] & & {K=[\Gg,B]}\ar@{-}[dl]^{C_4}\\
  & {K\times K=[K,D]}
}
\]

Let us define 
\[p=\PROJ adad,\qquad q=\PROJ ad,\qquad r=\PROJ a,\]
and consider the subset
$\Omega=\setsuch{Zp,Zq,Zr}{Z\in\{\DIAG,\PROJ\}^*}$ of $\Gg_*$, which
will serve as a basis for the quotients of successive terms of the
automorphism tower of $\Gg$. Define a rank on $\Omega$ by
\[\rank p=1,\qquad\rank q=2,\qquad\rank r=3\]
and inductively
\[\rank(\PROJ x)=2\rank x,\qquad\rank(\DIAG x)=2\rank x-1.\]
Then $\langle\Omega\cup S\rangle=\Gg_*$ is the layered closure of
$\Gg$.  Extend the rank function to $\Omega\cup S$ by $\rank(s)=0$ for
all $s\in S$, and define the following subgroups of $\Gg_*$:
\[A_n=\langle x\in\Omega\cup S:\rank x\le n\rangle.\]
The elements of $\Omega$ of low rank are as given below:
\[\begin{array}{r|ccccccccc}
  \text{rank} & \multispan{2}\text{\qquad elements}\\ \hline
  1 &      &       & p & \DIAG p & \DIAG\DIAG p & \DIAG\DIAG\DIAG p & \DIAG\DIAG\DIAG\DIAG p & \cdots & \DIAG^np\\
  2 &      & q     &   & \PROJ p    & \PROJ\DIAG p & \PROJ\DIAG\DIAG p & \PROJ\DIAG\DIAG\DIAG p & \cdots & \PROJ\DIAG^np\\
  3 & r    & \DIAG q  &   &      & \DIAG\PROJ p & \DIAG\PROJ\DIAG p & \DIAG\PROJ\DIAG\DIAG p & \cdots & \DIAG\PROJ\DIAG^np\\
  4 &      & \PROJ q  &   &      & \PROJ\PROJ p & \PROJ\PROJ\DIAG p & \PROJ\PROJ\DIAG\DIAG p & \cdots & \PROJ\PROJ\DIAG^np\\
  5 & \DIAG r & \DIAG\DIAG q &   &      &       & \DIAG\DIAG\PROJ p & \DIAG\DIAG\PROJ\DIAG p & \cdots & \DIAG\DIAG\PROJ\DIAG^np\\
  6 & \PROJ r & \PROJ\DIAG q &   &      &       & \PROJ\DIAG\PROJ p & \PROJ\DIAG\PROJ\DIAG p & \cdots & \PROJ\DIAG\PROJ\DIAG^np\\
  7 &      & \DIAG\PROJ q &   &      &       & \DIAG\PROJ\PROJ p & \DIAG\PROJ\PROJ\DIAG p & \cdots & \DIAG\PROJ\PROJ\DIAG^np\\
  8 &      & \PROJ\PROJ q &   &      &       & \PROJ\PROJ\PROJ p & \PROJ\PROJ\PROJ\DIAG p & \cdots & \PROJ\PROJ\PROJ\DIAG^np\\
  9 & \DIAG\DIAG r & \DIAG\DIAG\DIAG q &   &      &       &        & \DIAG\DIAG\DIAG\PROJ p & \cdots & \DIAG\DIAG\DIAG\PROJ\DIAG^np\\
  10& \PROJ\DIAG r & \PROJ\DIAG\DIAG q &   &      &       &        & \PROJ\DIAG\DIAG\PROJ p & \cdots & \PROJ\DIAG\DIAG\PROJ\DIAG^np\\
\end{array}\]
The next rows follow the same pattern: elements of rank $n$ are
determined by writing backwards $n-1$ in base $2$, using the symbols
$\DIAG=0$ and $\PROJ=1$.

\begin{thm}\label{thm:grigmain}
  We have for all $i\in\N$
  \[\aut^i(\Gg) = \langle a,b,c,d,x\in\Omega\text{ with }\rank x\le i\rangle.\]
  In particular, $\aut^{3\cdot2^{n-1}}(\Gg)$ contains $v*\Gg$ for all
  $v\in X^n$, since $\rank(\PROJ^na)=3\cdot2^{n-1}$.
  
  For all $i\in\N$, we have $\out(\aut^i(\Gg))=(\Z/2)^\omega$
  generated by the $x\in\Omega$ of rank $i+1$.
  
  We have $\aut^\alpha(\Gg)=\langle\Omega\rangle=\Gg_*$ for all
  $\alpha\ge\omega$.
\end{thm}

\begin{lem}\label{lem:grig-char}
  For all $s\in S$, the subgroups $\langle s\rangle^\Gg$ are all
  characteristic and distinct.
\end{lem}
\begin{proof}
  The level stabilizer $\stab_\Gg(n)$ is characteristic in $\Gg$ for
  all $n$, by Theorem~\ref{thm:ln}. The subgroup $\langle
  s\rangle^\Gg$ contains $\stab_\Gg(4)$, by direct
  computation~\cite{bartholdi-g:parabolic}*{Lemma~4.1}. To check that
  $\langle s\rangle^\Gg$ is normal in $\aut(\Gg)$, it is therefore
  sufficient to check that $\langle s\rangle^\Gg/\stab_\Gg(4)$ is
  normal in $\aut(\Gg)/\stab_\Gg(4)$.
  
  Now $\aut(\Gg)<W$ by Theorem~\ref{thm:ln}, so all computations may
  be performed in $\overline W=W/\stab_W(4)=\wr^4 C_2$, a $2$-Sylow
  subgroup of $\sym{2^4}$. The image $\overline\Gg$ of $\Gg$ has order
  $2^{12}$, and the image of $\aut(\Gg)$ lies in $\overline
  N=\norm_{\overline W}(\overline\Gg)$, a group of order $2^{13}$. It
  is then routine to check that $\langle s\rangle^\Gg/\stab_\Gg(4)$ is
  normal in $\overline N$.

  These computations were performed using the computer software
  \textsc{Gap}~\cite{gap4:manual}, by expressing $\overline W$ as a
  permutation group on $16$ points.\footnote{Other ad hoc arguments
    would be preferable for manual computation; for instance, list all
    quotients of $\Gg$ of order $8$ and $16$ to check that $D$ is
    characteristic, and that if $B,C$ are not characteristic, then
    there is an automorphism exchanging them. This last possibility is
    eliminated by computing the quotients $B/\stab_\Gg(3)=D_8\times
    C_2$ and $C/\stab_\Gg(3)=C_2^4$.}
\end{proof}

Define a norm on $\Gg_*$ by setting
\[\|g\|=\min\setsuch{n+\rank(\omega_1)+\dots+\rank(\omega_n)}{g=\omega_1\dots\omega_n,\,\omega_i\in\Omega\cup S}.\]
Note that the restriction of this norm to $\Gg$ coincides with the
norm given by~\eqref{eq:metric}, but we shall not need this fact.  We
have a contraction property on $\Gg_*$ as in $\Gg$:
\begin{lem}\label{lem:grig-contr}
  For all $g\in \Gg_*$ and $x\in X$ we have $|g@x\|\le\frac12(\|g\|+1).$
\end{lem}

Define next a norm on the set of group homomorphisms $\{\phi:\Gg\to \Gg_*\}$ by
\[|||\phi|||=\sum_{s\in S}\|\phi(s)\|.\]

\begin{lem}\label{lem:grig-g*}
  $\aut^\alpha(\Gg)\cong\norm^\alpha_{\Gg_*}(\Gg)$ for all ordinals
  $\alpha\le\omega$.
\end{lem}
\begin{proof}
  In~\cite{lavreniuk-n:rigidity}*{Theorem~8.5}, the authors show that
  $\Gg$ is saturated, with the subgroups $(\mho_2)^n(\Gg)$ fixing
  $X^n$ and acting transitively on each subtree at level $n$.
  Furthermore, $\Gg$ is weakly branch. By Lemma~\ref{lem:wbsat}, the
  groups $\aut^\alpha(\Gg)$ are also saturated and weakly branch, so
  Theorem~\ref{thm:ln} applies, and
  $\aut^{\alpha+1}(\Gg)\le\norm_W(\aut^\alpha(\Gg))$ for all
  ordinals $\alpha$.
  
  We proceed by induction on $\alpha$. Take $\phi\in\aut^\alpha(\Gg)$.
  First, we may suppose $\phi\in W$ by the above. We construct
  inductively elements $\phi_n\in W$ for all $n\in\N$ by setting
  $\phi_0=\phi$ and decomposing
  $\phi_n=(\lambda_n,\rho_n)a^{\varepsilon_n}$. We then note that
  $\lambda_n,\rho_n:\Gg\to \Gg_*$ since $\stab_\Gg(1)$ is subdirect in
  $\Gg\times \Gg$. We set $\phi_{n+1}$ to be that among
  $\lambda_n,\rho_n$ of smallest $|||$-norm, and proceed.
  
  We may assume $\varepsilon_n=0$, up to multiplying by an element on
  $\Gg_*$. Suppose $\rho_n$ is of minimal norm; then
  \begin{align*}
    |||\phi_{n+1}|||&=\|\phi_{n+1}(a)\|+\|\phi_{n+1}(b)\|+\|\phi_{n+1}(c)\|+\|\phi_{n+1}(d)\|\\
    &=\|\phi_n(axa)@\7\|+\|\phi_n(d)@\7\|+\|\phi_n(b)@\7\|+\|\phi_n(c)@\7\|\\
    &\le\frac12\left(2\|\phi_na\|+\|\phi_nx\|+\|\phi_nb\|+\|\phi_nc\|+\|\phi_nd\|+1\right)\\
    &<|||\phi_n|||
  \end{align*}
  for any $x\in\{b,c\}$; a similar equation holds if $\lambda_n$ is of
  minimal norm. Therefore $|||\phi_n|||=4$ for $n$ large enough, and
  by Lemma~\ref{lem:grig-char} we have $\phi_n=1$.
  
  Now for $m=n-1,n-2,\dots,0$ we have
  $[\PROJ \phi_m,a]=(\lambda_m^{-1}\rho_n,\rho_m\lambda_n^{-1})\in\aut^\alpha(\Gg)\subseteq
  \Gg_*$, and since one of $\lambda_m,\rho_m$ is $\phi_{m+1}$ and
  therefore belongs to $\Gg_*$, we obtain that both of them belong to
  $\Gg_*$, and therefore $\phi_m$ belongs to $\Gg_*$. We have shown
  $\phi_0\in \Gg_*$.
\end{proof}

\begin{proof}[Proof of Theorem~\ref{thm:grigmain}]
  By Lemma~\ref{lem:grig-g*}, the automorphism tower of $\Gg$ is
  contained in $\Gg_*$. We prove the theorem by showing that the
  generators of $\out(A_n)$ are precisely the elements of rank $n+1$; 
  for this purpose, we show that
  \begin{enumerate}
  \item elements of rank $n+1$ normalize $A_n$, and
  \item no element of rank $n+2$ or greater normalizes $A_n$.
  \end{enumerate}
  It then follows that $A_n=\norm_{\Gg_*}^n(\Gg)$.
  
  The following inductive computations show that
  $\rank[x,y]\le\max\{\rank x,\rank y\}-1$ for $x\in\Omega$ and $y\in
  \Omega\cup \Gg$:
  \begin{align*} 
    \rank[\PROJ x,\PROJ y]&=\rank\PROJ [x,y]=2\rank[x,y]\le2(\max\{\rank
    x,\rank y\}-1)\\
    &\le\max\{\rank\PROJ x,\rank\PROJ y\}-1;\\
    \rank[\DIAG x,\PROJ y]&=\rank\PROJ [x,y]=2\rank[x,y]\le2(\max\{\rank
    x,\rank y\}-1)\\
    &\le\max\{\rank\DIAG x,\rank\PROJ y\}-1;\\
    \rank[\DIAG x,\DIAG y]&=\rank\DIAG [x,y]\\
    &=2\rank[x,y]-1\le2(\max\{\rank x,\rank y\}-1)-1\\
    &\le\max\{\rank\DIAG x,\rank\DIAG y\}-1;\\
    \rank[\PROJ x,a]&=\rank\DIAG x=\rank(\PROJ x)-1;\\
    \rank[\DIAG x,a]&=\rank1=0;\\
    \rank[\PROJ x,b]&=\rank\PROJ [x,a]=2\rank[x,a]\le2(\rank x-1)\le\rank(\PROJ x)-1;\\
    \rank[\DIAG x,b]&=\rank(\PROJ [x,a])(1,[x,c])\le\max\{\rank[x,a],\rank[x,c]\}\\
    &\le\rank(\DIAG x)-1;\\
    [p,a]&=(ac)^8\in \Gg,\text{ and more generally }[p,s]\in \Gg\text{ for }
    s\in S;\\
    \rank[\PROJ x,r]&=\rank\PROJ [x,adad]\le\rank(\PROJ x)-1;\\
    [q,s]&\in \Gg\langle p\rangle\text{ for }s\in S;\\
    [r,s]&\in \Gg\langle q\rangle\text{ for }s\in S;\\
    [p,q]&=1,\quad [p,r]=1,\quad [q,r]=p.
  \end{align*}
  
  Finally, the tower of $\Gg$ stops at level $\omega$ by
  Theorem~\ref{thm:stops@omega}.
%
\end{proof}

We note by direct computation that the ``branching subgroup'' $K$
satisfies the property that $\times^i(K)$ is characteristic in
$\aut^{2^i}(\Gg)$. It then follows that $\Gg\cap \Gg^\alpha$ has
finite index in $\Gg$ for all $\alpha\in\aut^i(G)$.

More generally, we have shown for $G\in\{\Gg,\GS\}$ that
$\aut^{i+1}(G)/\aut^i(G)$ is solvable (actually, of degree $\leq 2$)
for all $i\in\N$. Therefore, there exists $k\in\N$ such that
$(\aut^i(G))^{(k)}\leq G$, where $H^{(k)}$ denotes the $k$th term in
the derived series of $H$. Since $G$ is just-infinite and not
solvable, it follows that $G\cap G^\alpha$ has finite index in both
$G$ and $G^\alpha$, for all $\alpha \in\aut^\omega(G)$. Hence,
$\aut^\omega(G)$ is a subgroup of the commensurator of $G$ in $W$.

\section{Thanks}
The authors express their gratitude to the referee who carefully read
this manuscript, and pointed out the relevance of Fodor's lemma to
Proposition~\ref{prop:towercountable}.
\end{document}